\documentclass[12pt,a4paper]{article}
\usepackage{graphicx}
\usepackage{amsmath,amsthm,amssymb,amscd,mathtools,amsfonts}
\usepackage[noblocks]{authblk}
\usepackage{mathrsfs}
\usepackage{titlesec}
\usepackage{appendix}
\usepackage{enumerate}
\usepackage{epsfig}
\usepackage[square, comma, sort&compress, numbers]{natbib}
\usepackage{epstopdf}
\usepackage{bm}
\usepackage[T1]{fontenc}
\usepackage[utf8]{inputenc}
\usepackage{authblk}
\usepackage{tabularx}
\newcommand{\RNum}[1]{\uppercase\expandafter{\romannumeral #1\relax}}
\begin{document}
\title{A New Class of Multivariate Elliptically Contoured Distributions with Inconsistency Property}
\author[a]{Yeshunying Wang}
\author[a]{Chuancun Yin}
\affil[a]{School of Statistics, Qufu Normal University \authorcr Shandong 273165, China \authorcr e-mail: ccyin@qfnu.edu.cn}
\renewcommand*{\Affilfont}{\small\it}

\maketitle
\begin{abstract}
 We introduce a new class of multivariate elliptically symmetric distributions including elliptically symmetric logistic distributions and Kotz type distributions. We investigate the various probabilistic properties including marginal distributions, conditional distributions, linear transformations, characteristic functions and dependence measure in the perspective of the inconsistency property. In addition, we provide a real data example to show that the new distributions  have reasonable flexibility.
\end{abstract}
\newpage
\section{Introduction}
\par
The multivariate normal distribution has been widely used in theory and practice because of its tractable statistical features. However, the light tail of the normal distribution can not fit some practical situation well. The elliptically contoured distributions (elliptical distributions), a new family of distributions with similar convenient properties, overcomes the shortcomings of the normal distributions. An $n$-dimension random vector \textbf{X} is said to have a multivariate elliptical distribution, written as
  $\textbf{X}\sim Ell_{n}(\bm{\mu},\mathbf{\Sigma},\phi)$ if its characteristic function can be expressed as $\psi_{\textbf{X}}(\textbf{t})=\exp(i\textbf{t}^{\emph{T}}\bm{\mu})\phi(\textbf{t}^{\emph{T}}\mathbf{\Sigma} \textbf{t})$, where $\bm{\mu}$ is an $n$-dimension column vector, $\mathbf{\Sigma}$ is an $n\times n$ positive semi-definite matrix, $\phi(\cdot)$ is called characteristic generator. If $\textbf{X}$ has a probability density function (pdf) $f(\textbf{x})$, then
  $$f(\textbf{x})=\frac{C_{n}}{\sqrt{|\mathbf{\Sigma}|}}
                      g_{n}\left((\textbf{x}-\bm{\mu})^{\emph{T}}\mathbf{\Sigma}^{-1}(\textbf{x}-\bm{\mu})\right),$$
where $C_{n}$ is the normalizing constant and $g_{n}(\cdot)$ is called density generator (d.g.).
The stochastic representation of \textbf{X} is given by
$$\textbf{X}=\bm{\mu}+R\textbf{A}^{\emph{T}}\textbf{U}^{(n)},\eqno(1.1)$$ where \textbf{A} is a square matrix such that $\textbf{A}^{\emph{T}}\textbf{A}=\mathbf{\Sigma}$, $\textbf{U}^{(n)}$ is uniformly distributed on the unit sphere surface in $\mathbb{R}^{n}$ , $R\geq 0$ is independent of $\textbf{U}^{(n)}$ and has the pdf given by $$f_{R}(v)=\frac{1}{\int_{0}^{\infty}t^{n-1}g(t^{2})dt}v^{n-1}g(v^{2}),~v\geq0.\eqno(1.2)$$
\par
Many members of the elliptical distributions such as the multivariate normal distributions and student-$t$ distributions, have been systematic studied. See the books and papers of Cambanis et al. (1981), Fang et al. (1990), Kotz and Ostrovskii (1994), Liang and Bentler (1998), Nadarajah (2003).
Nevertheless, research work on the multivariate symmetric logistic distribution is far less than other members. The elliptically symmetric logistic distribution with density
$$f(\textbf{x})=\frac{\Gamma(\frac{n}{2})|\mathbf{\Sigma}|^{-\frac{1}{2}}}
            {\pi^{\frac{n}{2}}\int_{0}^{\infty}u^{\frac{n}{2}}\frac{\exp(-u)}{(1+\exp(-u))^{2}}du}
\frac{\exp(-(\textbf{x}-\bm{\mu})^{\emph{T}}\mathbf{\Sigma}^{-1}(\textbf{x}-\bm{\mu}))}
     {[1+\exp(-(\textbf{x}-\bm{\mu})^{\emph{T}}\mathbf{\Sigma}^{-1}(\textbf{x}-\bm{\mu}))]^2},~\textbf{x}\in \mathbb{R}^{n},
$$was introduced by Jensen (1985) and has been studied by Fang et al. (1990), Kano (1994), Yin and Sha (2018). Several applications of multivariate symmetric logistic distribution in risk management, quantitative finance and actuarial science can be found in various literatures such as Landsman and Valdez (2003), Landsman et al. (2016a, 2016b, 2018).
\par
The paper will define a new class of elliptical distributions including the Kotz type distributions and the logistic distributions, give the value of the normalizing constant and study the marginal distributions, conditional distributions, linear transformations, characteristic functions and local dependence functions in perspective of its inconsistency property.
\par
The rest of the paper is organized as follows. In Section 2, we introduce the definition of a new class of  multivariate elliptically symmetric distributions which include elliptically symmetric logistic distributions and Kotz type distributions and discuss the expression of the normalizing constant. In Sections 3-7 we study the probabilistic properties of the new class of  elliptically distribution including marginal distributions, conditional distributions, linear transformations, characteristic functions in perspective of its inconsistency property. In addition, we give the expression of its local dependence function.
In Section 8, we give the data analysis of the  new class of  elliptically distribution and we conclude in Section 9.
\section{Preliminary}
We now give the definition of a new class of  multivariate elliptically symmetric distributions which include elliptically symmetric logistic distributions and Kotz type distributions.
\par
\begin{flushleft}
\textbf{Definition 2.1}
  The $n$-dimensional random vector $\textbf{X}$ is said to have a generalized elliptical logistic (GL) distribution  with parameter $\bm{\mu}$ $($$n$-dimensional vector$)$ and  $\mathbf{\Sigma}$ $($$n\times n$ matrix with $\mathbf{\Sigma}$$>$0$)$ if its pdf and density generator have the forms
   $$
   f(\textbf{x})= C_{n}\mid\mathbf{\Sigma}\mid^{-\frac{1}{2}}g\left((\textbf{x}-\bm{\mu})^{\emph{T}}
   \mathbf{\Sigma}^{-1}(\textbf{x}-\bm{\mu})\right),~
   \textbf{x}\in \mathbb{R}^{n},\eqno(2.1)$$
$$ g(t)=\frac{t^{N-1}\exp(-at^{s_{1}})}{\left(1+\exp(-bt^{s_{2}})\right)^{2r}},~t>0,\eqno(2.2)$$ respectively, where $2N+n>2$, $a,~b,~ s_{1},~ s_{2}>0$, $r\geq0$ are constants. The normalizing constant $C_{n}$ will be discussed in Section 2.2.
\end{flushleft}
\subsection{Special cases}
\begin{enumerate}[1)]
  \item \textbf{Generalized logistic distribution} (Yin and Sha (2018))
  \\Setting $N=1,~s_{1}=s_{2}=1$ in $(2.2)$, we get $$g(t)=\frac{\exp(-at)}{(1+\exp(-bt))^{2r}},\eqno(2.3)$$which is the density generator put forward by Yin and Sha (2018).
  \item \textbf{Multivariate normal distribution}\\
   Setting $N=1,~a=\frac{1}{2},~s_{1}=1,~r=0$ in $(2.2)$, we get $g(t)=\exp(-\frac{1}{2}t)$, which is the density generator of the multivariate normal distribution.
   \item\textbf{Multivariate exponential power (Epo) distribution} (Landsman and Valdez (2003))\\
   For $N=1,~r=0$, $(2.2)$ is the density generator of the multivariate exponential power distribution whose d.g. is usually written as$$g(t)=\exp(-at^{s_{1}}),~a,~s_{1}>0.$$
   If $s_{1}=\frac{1}{2}$ and $a=\sqrt{2}$, we have the d.g. of the double exponential or Laplace distribution defined as
   $$g(t)=\exp(-\sqrt{2t}).$$
 \item \textbf{Kotz type (Ko) distribution} (Fang et al. (1990))\\
   For $r=0$, (2.2) is the density generator of the symmetric Kotz type distribution whose d.g. is usually written as$$g(t)=t^{N-1}\exp(-at^{s_{1}}),~ a, ~s_{1}>0, ~2N+n>2.\eqno(2.4)$$
 When $s_{1}=1$, $(2.4)$ is the density generator of the original Kotz distribution whose d.g. is written as
  $$g(t)=t^{N-1}\exp(-at), ~a>0, ~2N+n>2.$$
   \item \textbf{Elliptically symmetric logistic (Lo) distribution} (Fang et al. (1990))\\
 Setting $N=1, ~a=b=r=1,~ s_{1}=s_{2}=1$, $(2.2)$ is the d.g. of the $n$-dimensional elliptically symmetric logistic distribution , written as $$g_{n}(t)=\frac{\exp(-t)}{(1+\exp(-t))^{2}}.\eqno(2.5)$$
  \item \textbf{Generalized logistic type \RNum{1} (GL\RNum{1}) distribution} (Arashi and Nadarajah (2016))\\
  Setting $N=1,~a=b=1,~s_{1}=s_{2}=1$ in $(2.2)$ gives the density generator of the generalized logistic type \RNum{1} distribution
  written as $$g(t)=\frac{\exp(-t)}{(1+\exp(-t))^{2r}}.\eqno(2.6)$$
  \item \textbf{Generalized logistic type \RNum{3} (GL\RNum{3}) distribution} (Arashi and Nadarajah (2016))\\
 Setting $N=1,~b=1,~s_{1}=s_{2}=1,~r=a$ in $(2.2)$ gives the density generator of the generalized logistic type \RNum{3} distribution
  written as $$g(t)=\frac{\exp(-at)}{(1+\exp(-t))^{2a}}.\eqno(2.7)$$
  \item \textbf{Generalized logistic type \RNum{4} (GL\RNum{4}) distribution} (Arashi and Nadarajah (2016))\\
  For $N=1,~s_{1}=s_{2}=1$, $r=\frac{p+a}{2}$, where $p>0$, (2.2) is the density generator of the generalized logistic type \RNum{4} distribution written as
  $$g(t)=\frac{\exp(-at)}{(1+\exp(-t))^{p+a}}.\eqno(2.8)$$
\end{enumerate}
Setting $s_{1}=s_{2}=s$ in $(2.2)$, we obtain $$g_{n}(t)=\frac{t^{N-1}\exp(-at^{s})}{(1+\exp(-bt^{s}))^{2r}}.\eqno(2.9)$$
\par
For the sake of simplicity, we discuss probabilistic properties of GL distributions with density generators defined as (2.9) in following sections.
\subsection{Normalizing constant}
\par
To calculate the normalizing constant defined in $(2.1)$, we introduce the Hurwitz-Lerch zeta function. The Hurwitz-Lerch zeta function and its integral representation  are respectively defined as
$$\Phi(z,s,a)=\sum_{n=0}^{\infty}\frac{z^{n}}{(n+a)^{s}}$$
$$(a\in \mathbb{C}\setminus\mathbb{Z}_{0}^{-},s\in C~when\mid z\mid<1;\Re(s)>1~when\mid z\mid=1 ),$$
$$\Phi(z,s,a)=\frac{1}{\Gamma(s)}\int_{0}^{\infty}\frac{t^{s-1}e^{-at}}{1-ze^{-t}}dt=\frac{1}{\Gamma(s)}\int_{0}^{\infty}\frac{t^{s-1}e^{-(a-1)t}}{e^{t}-z}dt$$
$$(\Re(s)>0,\Re(a)>0~when\mid z\mid\leq1(z\neq1);\Re(s)>1~when~z=1).$$
Various generalization and extensions of the Hurwitz-Lerch zeta function $\Phi(z,s,a)$ have been studied by various researchers. The following expression of generalized Hurwitz-Lerch zeta function which will be used in this paper is defined by (cf. Lin et al. (2006))
$$\Phi_{v}^{*}(z,s,a)=\frac{1}{\Gamma(v)}\sum_{n=0}^{\infty}\frac{\Gamma(v+n)}{n!}\frac{z^{n}}{(n+a)^{s}}$$
$$(v\in\mathbb{C},a\in \mathbb{C}\setminus\mathbb{Z}_{0}^{-},s\in\mathbb{C}~when \mid z\mid<1;\Re(s-v)>1~when\mid z\mid=1),$$
$$\Phi_{v}^{*}(z,s,a)=\frac{1}{\Gamma(s)}\int_{0}^{\infty}\frac{t^{s-1}e^{-at}}{(1-ze^{-t})^{v}}dt=\frac{1}{\Gamma(s)}\int_{0}^{\infty}\frac{t^{s-1}e^{-(a-v)t}}{(e^{t}-z)^{v}}dt$$
$$(\Re(s)>0,\Re(a)>0~when\mid z\mid\leq1(z\neq1);\Re(s)>1~when~z=1).$$
If $v=0$, $$\Phi_{v}^{*}(z,s,a)=\Phi_{0}^{*}(z,s,a)=\frac{1}{a^{s}},$$
thus we denote $\Phi_{0}^{*}(z,s,a)$ by $\Phi_{0}^{*}(s,a)$.
\par
Pointed out by Yin and Sha (2018), the normalizing constant of elliptical symmetric logistic distribution suggested by Landsman and Valdez (2003), has no meaning when $n$=1 and $n$=2. Thus, we calculate the normalizing constant defined in (2.1) by the generalized Hurwitz-Lerch zeta function. The method was similarly used in Yin and Sha (2018).
\par
\begin{flushleft}
\textbf{Theorem 2.1} Letting $\textbf{X}\sim GL_{n}(\bm{\mu}, \mathbf{\Sigma}, g_{n})$ where $g_{n}$ is defined as $(2.9)$, then the normalizing constant defined in $(2.1)$ can be expressed as
\end{flushleft}
$$ C_{n}=c_{n}^{*}(N,b,s) \left[\Phi_{2r}^{*}(-1,\frac{1}{s}(N+\frac{n}{2}-1),\frac{a}{b})\right]^{-1},$$
        where $$~c_{n}^{*}(N,b,s)=\frac{\Gamma(\frac{n}{2})}{\Gamma(\frac{1}{s}(N+\frac{n}{2}-1))}
        \frac{b^{\frac{1}{s}(N+\frac{n}{2}-1)}s}{\pi^{\frac{n}{2}}},$$ and $\Phi_{2r}^{*}$ is the generalized Hurwitz-Lerch zeta function.
\par
\textbf{Proof.} Since $$\int_{-\infty}^{\infty}\int_{-\infty}^{\infty}\cdots\int_{-\infty}^{\infty}f(\mathbf{x})d\mathbf{x}=1,$$ where $f(\mathbf{x})$ is defined in $(2.1)$ and transformation from the rectangular to polar coordinates. We have
\begin{equation}
\begin{split}
C_{n}&=\frac{\Gamma(\frac{n}{2})}{\pi^{\frac{n}{2}}}\
             \left[\int_{0}^{\infty}x^{\frac{n}{2}-1}g_{n}(x)dx\right]^{-1}\nonumber\\
      &=\frac{\Gamma(\frac{n}{2})}{\pi^{\frac{n}{2}}}\
              \left[\int_{0}^{\infty}\frac{x^{N+\frac{n}{2}-2}e^{-ax^{s}}}{(1+e^{-bx^{s}})^{2r}}dx\right]^{-1}\\
      &=\frac{\Gamma(\frac{n}{2})}{\pi^{\frac{n}{2}}}
        \frac{b^{\frac{1}{s}(N+\frac{n}{2}-1)}s}{\Gamma(\frac{1}{s}(N+\frac{n}{2}-1))}\left[\Phi_{2r}^{*}\left(-1,\frac{1}{s}(N+\frac{n}{2}-1),\frac{a}{b}\right)\right]^{-1}\\
      &=c_{n}^{*}(N,b,s) \left[\Phi_{2r}^{*}\left(-1,\frac{1}{s}(N+\frac{n}{2}-1),\frac{a}{b}\right)\right]^{-1},
\end{split}
\end{equation}
where $$c_{n}^{*}(N,b,s)=\frac{\Gamma(\frac{n}{2})}{\Gamma(\frac{1}{s}(N+\frac{n}{2}-1))}
        \frac{b^{\frac{1}{s}(N+\frac{n}{2}-1)}s}{\pi^{\frac{n}{2}}}.$$
\begin{flushleft}
\textbf{Corollary 2.1.1}
\begin{enumerate}[1)]
  \item Supposing $\textbf{X}\sim Ko_{n}(\bm{\mu},\mathbf{\Sigma}, g_{n})$, where $g_{n}$ is defined as $(2.4)$, then the normalizing constant defined in $(2.1)$ can be expressed as $$ C_{n}=\frac{s_{1}\Gamma(\frac{n}{2})}{\pi^{\frac{n}{2}}\Gamma(\frac{1}{s_{1}}(N+\frac{n}{2}-1))}
        \left[\Phi_{0}^{*}\left(\frac{1}{s_{1}}(N+\frac{n}{2}-1), a\right)\right]^{-1}.$$
  \item Supposing $\textbf{X}\sim Epo_{n}(\bm{\mu},\mathbf{\Sigma}, g_{n})$, where $g_{n}$ is
  $$g_{n}(t)=\exp(-at^{s_{1}}),~a,~s_{1}>0,$$ then the normalizing constant defined in $(2.1)$ can be expressed as
  $$ C_{n}=
        \frac{s\Gamma(\frac{n}{2})}{\pi^{\frac{n}{2}}\Gamma(\frac{n}{2s_{1}})}\left[\Phi_{0}^{*}(\frac{n}{2s_{1}}, a)\right]^{-1}.$$
  \item Supposing $\textbf{X}\sim GL\RNum{1}_{n}(\bm{\mu},\mathbf{\Sigma}, g_{n})$, where $g_{n}$ is defined as $(2.9)$,
   then the normalizing constant defined in $(2.1)$ can be expressed as $[\pi^{\frac{n}{2}}\Phi_{2r}^{*}(-1,\frac{n}{2},1)]^{-1}$.
   \item Supposing $\textbf{X}\sim GL\RNum{3}_{n}(\bm{\mu},\mathbf{\Sigma}, g_{n})$, where $g_{n}$ is defined as $(2.7)$,
   then the normalizing constant defined in $(2.1)$ can be expressed as    $[\pi^{\frac{n}{2}}\Phi_{2a}^{*}(-1,\frac{n}{2},a)]^{-1}$.
    \item Supposing $\textbf{X}\sim GL\RNum{4}_{n}(\bm{\mu},\mathbf{\Sigma}, g_{n})$, where $g_{n}$ is defined as $(2.8)$,
   then the normalizing constant defined in $(2.1)$ can be expressed as
      $[\pi^{\frac{n}{2}}\Phi_{p+a}^{*}(-1,\frac{n}{2},a)]^{-1}$.
\end{enumerate}
\end{flushleft}
\begin{figure}[htbp]
\centering
\includegraphics[scale=0.3]{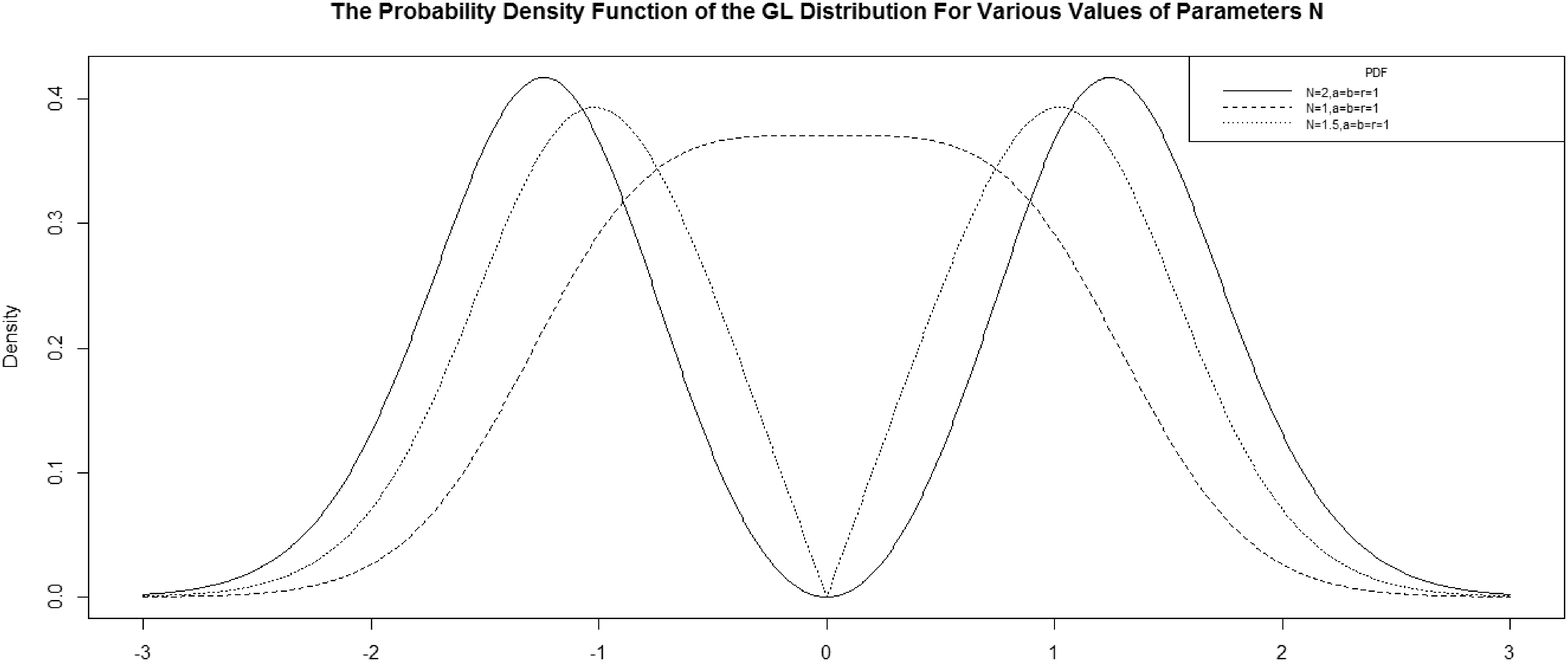}
\includegraphics[scale=0.3]{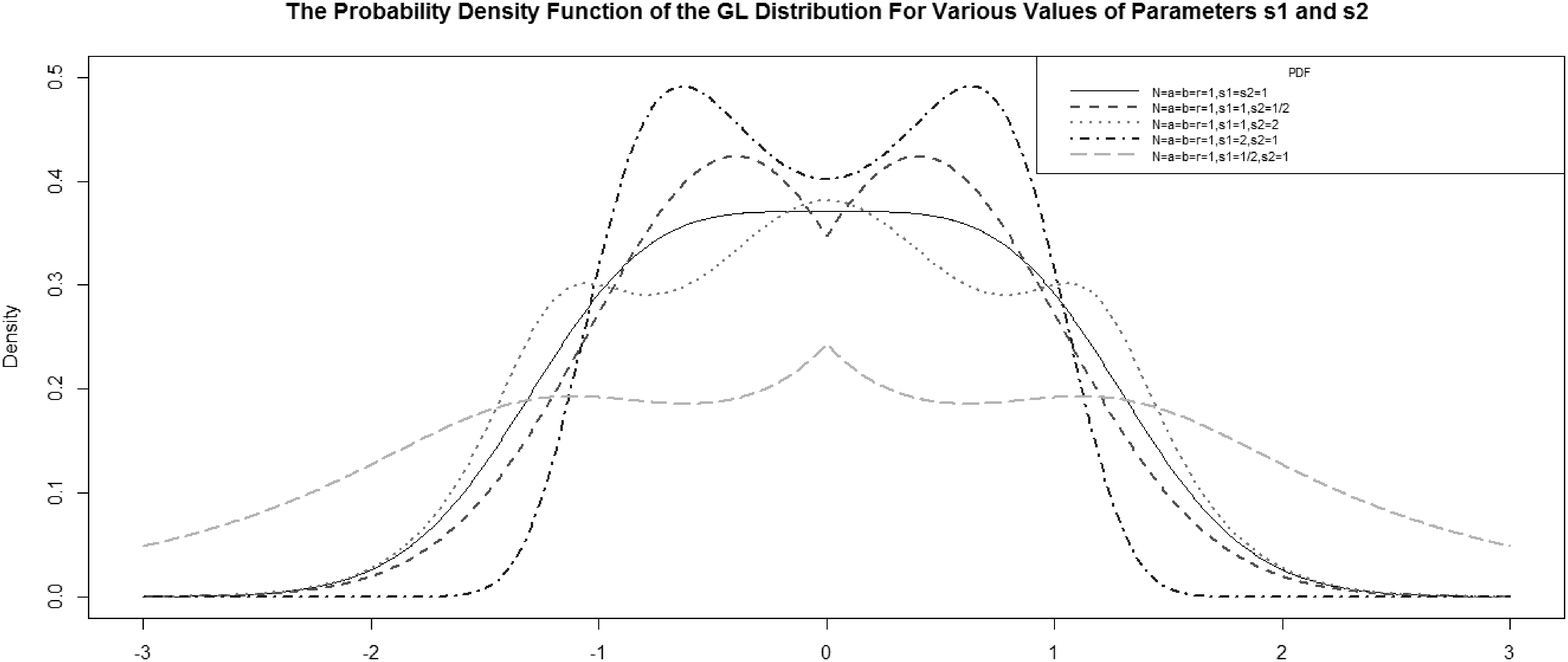}
\end{figure}
\begin{figure}[htbp]
\includegraphics[scale=0.3]{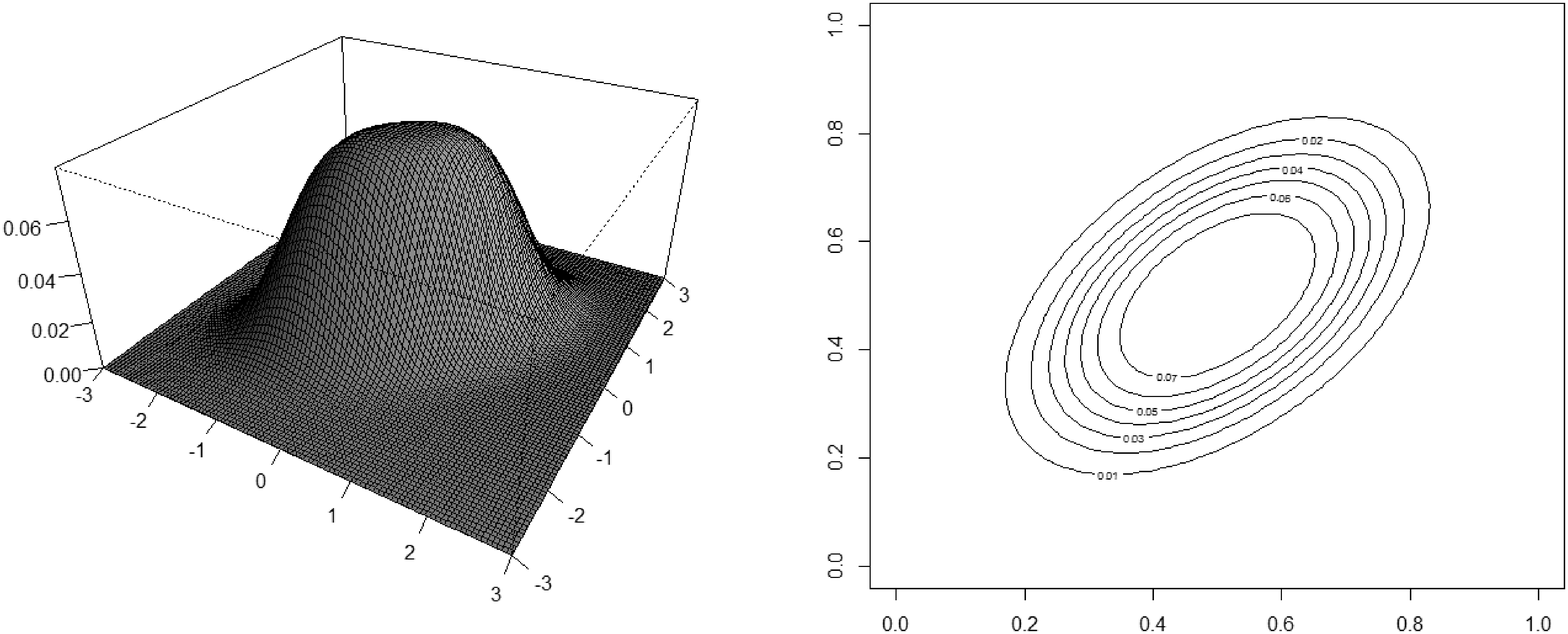}
\caption{$a=b=r=1,~N=1,~s_{1}=s_{2}=1,~\rho=0.5$.}
\includegraphics[scale=0.3]{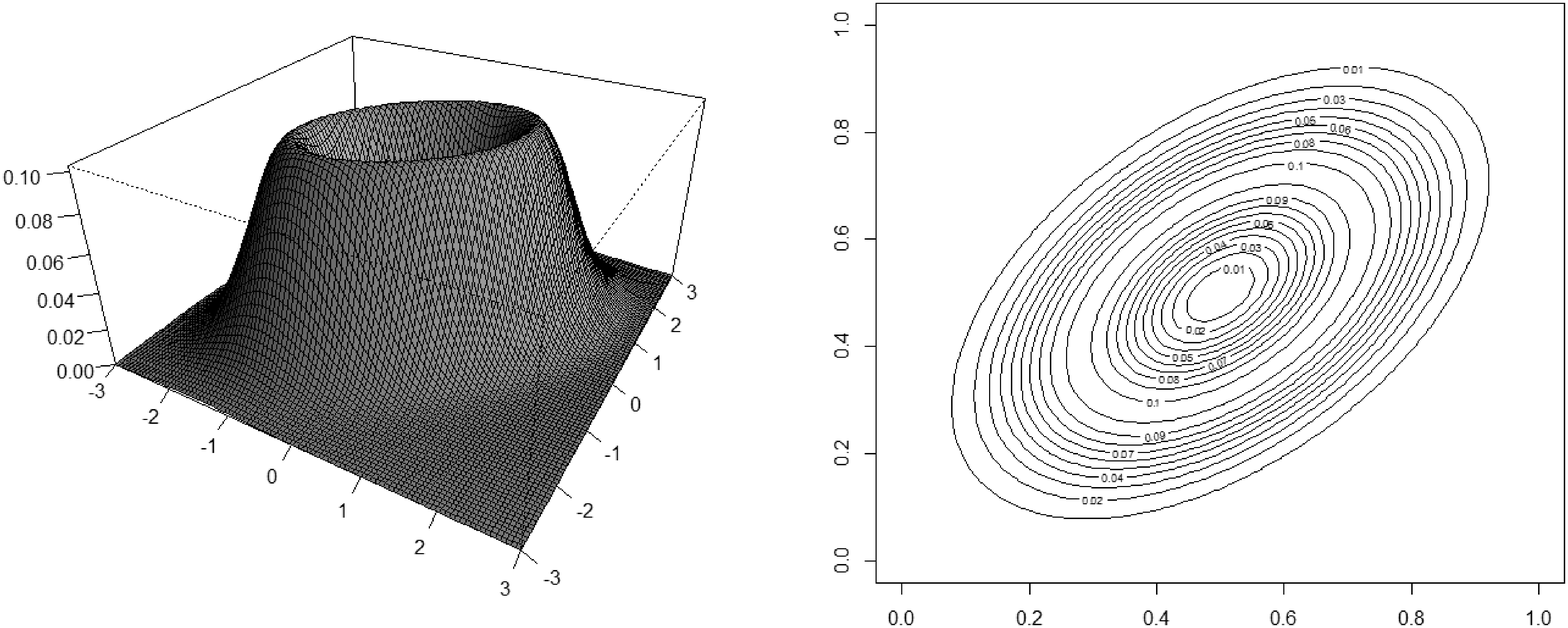}
\caption{$a=b=r=1,~ N=2,~ s_{1}=s_{2}=1,~\rho=0.5$.}
\includegraphics[scale=0.3]{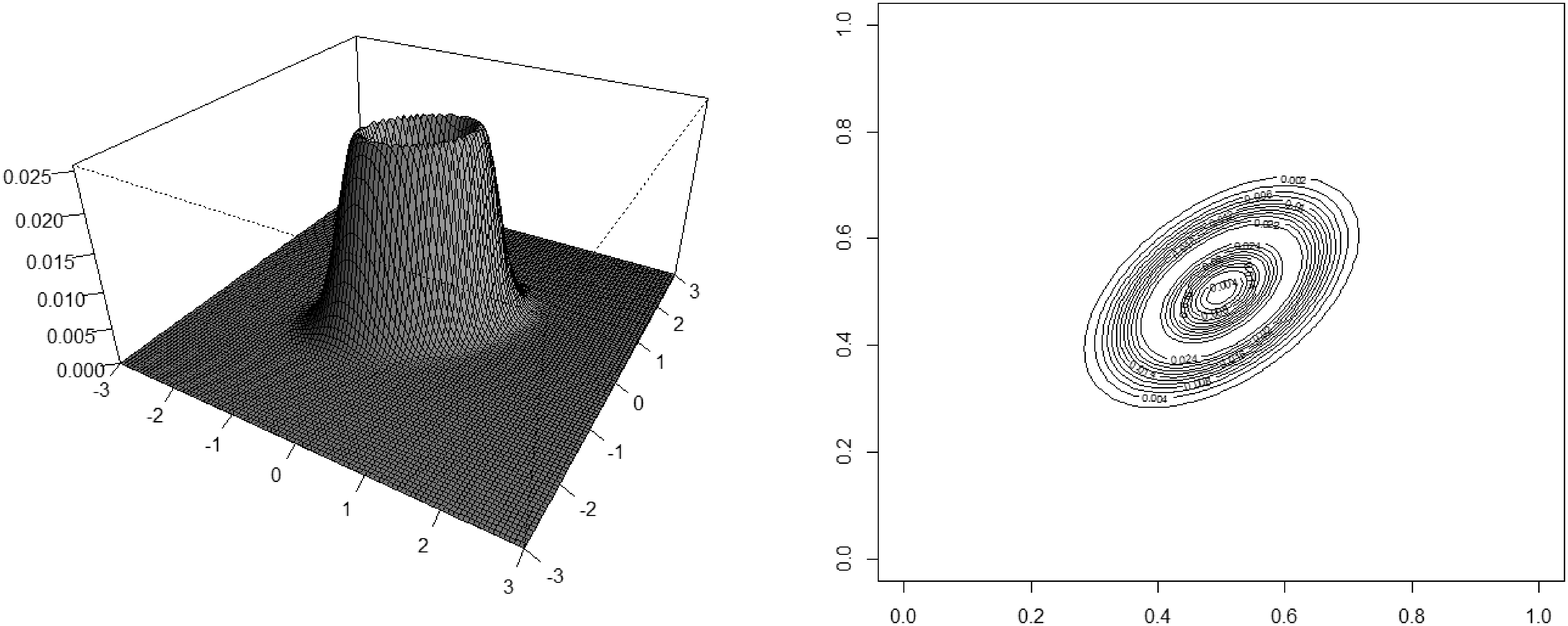}
\caption{$a=b=4,~ r=1, ~N=2,~ s_{1}=s_{2}=1,~ \rho=0.5$.}
\end{figure}
\par
Consider a family of density generators $$\{f(u|p)|p\in \mathbb{N}\},\eqno(3.1)$$ where $\mathbb{N}$ denotes the set of all positive integers. According to Kano $(1994)$, we will say that the family in $(3.1)$ possesses a consistency property if and only if $$\int_{-\infty}^{\infty}f\left(\sum_{j=1}^{p+1}x_{j}^{2}|p+1\right)dx_{p+1}=f\left(\sum_{j=1}^{p}x_{j}^{2}|p\right)\eqno(3.2)$$for any $p\in \mathbb{N}$ and almost all $(x_{1},\cdots,x_{p})\in \mathbb{R}^{p}$. We also say that the family is dimension coherent. However, the family of density generators defined in (2.2) does not satisfy the consistency property i.e. the marginal distributions of the $n$-dimension $(n>1)$ generalized elliptical bimodal logistic distribution --- which are elliptically contoured --- don't have (2.2) as their density generators. In the following sections we will discuss the inconsistency property of the GL distribution and its applications.
\section{Density generator}
\par
 Supposing $g_{1}(x)$, the d.g. of GL distributed random variable $\textbf{X}$, is defined as follows
 $$g_{1}(t)=\frac{t^{N-1}\exp(-at)}{(1+\exp(-bt))^{2r}}.$$
 Before we utilize (cf. Fang et al. (1990))
 $$g_{m}(u)=\int_{u}^{\infty}(\omega-u)^{\frac{n-m}{2}-1}g_{n}(\omega)d\omega,~g_{m}(u)=\int_{u}^{\infty}g_{m+2}(\omega)d\omega,\eqno(3.3)$$ to investigate $g_{2},~g_{3},~\cdots~,$ it is sufficient to verify that $g_{1}$ is a non-increasing function i.e. $g_{1}'(u)\leq0$. Without loss of generality, assume $a=b=1$, then
$$g_{1}'(u)=\frac{(N-1)u^{N-2}e^{-u}-t^{N-1}e^{-u}+(N-1)u^{N-2}e^{-2u}-t^{N-1}e^{-2u}+2ru^{N-1}e^{-2u}}{(1+e^{-u})^{2r+1}}.$$ $g_{1}'(u)\leq0$ if and only if
\begin{equation}
\begin{split}
&\frac{u^{N-2}e^{-u}}{(1+e^{-u})^{2r+1}}[(N-1)(u+e^{-u})-u(1+e^{-u})+2rue^{-u}]\leq0,\\ \nonumber
&N\leq\frac{u+e^{-u}+u+ue^{-u}-2rue^{-u}}{u+e^{-u}}.
\end{split}
\end{equation}
When $u=0$, $N\leq1$; $0<u<1$, $N\leq1$; $u\geq1$, $N\leq2$.
Above all, $g_{1}$ is a non-increasing function when $N\leq1$.
 Therefore, if $N\leq1$,
 $$\frac{t^{N-1}\exp(-at)}{(1+\exp(-bt))^{2r}}=\int_{t}^{\infty}g_{3}(u)du,~
 g_{2}(u)=\int_{t}^{\infty}(u-t)^{-\frac{1}{2}}g_{3}(u)du.$$
 Then, we have \begin{equation}
\begin{split}
g_{3}(u)=&-\frac{1}{(1+e^{-bu})^{2r+1}}\left[(N-1)u^{N-2}e^{-au}-au^{N-1}e^{-au}+(N-1)u^{N-2}e^{-(a+b)u}\right.\nonumber\\
&\left.-au^{N-1}e^{-(a+b)u}+2bru^{N-1}e^{-(a+b)u}\right],
\end{split}
\end{equation}
\begin{equation}
\begin{split}
g_{2}(u)=&ae^{-au}\sum_{j=0}^{N-1}\frac{(N-1)_{j}}{j!}\frac{\Gamma(j+\frac{1}{2})}{b^{j+\frac{1}{2}}}u^{N-1-j}\Phi_{2r+1}^{*}(-e^{-bu},j+\frac{1}{2},\frac{a}{b})\nonumber\\
         &-(N-1)e^{-au}\sum_{k=0}^{N-1}\frac{(N-2)_{k}}{k!}\frac{\Gamma(k+\frac{1}{2})}{b^{k+\frac{1}{2}}}u^{N-2-k}\Phi_{2r+1}^{*}(-e^{-bu},k+\frac{1}{2},\frac{a}{b})\\
         &-(N-1)e^{-(a+b)u}\sum_{l=0}^{N-1}\frac{(N-2)_{l}}{l!}\frac{\Gamma(l+\frac{1}{2})}{b^{l+\frac{1}{2}}}u^{N-2-l}\Phi_{2r+1}^{*}(-e^{-bu},l+\frac{1}{2},\frac{a+b}{b})\\
         &+ae^{-au}\sum_{v=0}^{N-1}\frac{(N-1)_{v}}{v!}\frac{\Gamma(v+\frac{1}{2})}{b^{v+\frac{1}{2}}}u^{N-1-v}\Phi_{2r+1}^{*}(-e^{-bu},v+\frac{1}{2},\frac{a+b}{b})\\
         &-2bre^{-(a+b)u}\sum_{q=0}^{N-1}\frac{(N-1)_{q}}{q!}\frac{\Gamma(q+\frac{1}{2})}{b^{q+\frac{1}{2}}}u^{N-1-q}\Phi_{2r+1}^{*}(-e^{-bu},q+\frac{1}{2},\frac{a+b}{b}),
\end{split}
\end{equation}
where $(x)_{n}=x(x-1)(x-2)\cdots(x-n+1)$.
In the same way, we can obtain $g_{5},~g_{7},~\cdots,$ if $g_{1}$ is a completely monotone function i.e. $(-1)^{n}g_{1}^{(n)}(t)\geq0$ for every $n\in \mathbb{N}_{0}$ and $t>0$, where $\mathbb{N}_{0}$ denotes the set of non-negative integers. After that we can obtain $g_{4},~g_{6},~\cdots$.
If $\textbf{X}$ has an elliptically symmetric logistic distribution i.e. $\textbf{X} \sim Lo_{n}(\bm{\mu},\mathbf{\Sigma},g_{n})$, where $g_{n}$ is defined as (2.5), we can obtain
$$g_{2}(t)=\sqrt{\pi}(\Phi_{3}^{*}(-e^{-t},\frac{1}{2},1)-\frac{e^{-t}}{\sqrt{2}}\Phi_{3}^{*}(-e^{-t},\frac{1}{2},2)),$$
$$g_{3}(t)=\frac{e^{-t}-e^{-2t}}{(1+e^{-t})^{3}}.$$
\par
In the same way, if $\textbf{X}=(X_{1},\cdots, X_{n})^{\emph{T}} \sim GL_{n}(\bm{\mu},\mathbf{\Sigma},g_{n})$, the d.g. of $X_{i}~(i=1,2,\cdots,n)$ differs from $g_{n}$ obviously.
We conclude the property in the following \\theorem.
\begin{flushleft}
\textbf{Theorem 3.1}~Letting $\textbf{X}=(\textbf{X}_{(m)}^{\emph{T}},\textbf{X}_{(n-m)}^{\emph{T}})^{\emph{T}} \sim GL_{n}(\bm{\mu},\mathbf{\Sigma},g_{n})$, where $g_{n}(t)$ is defined as $(2.2)$, $1\leq m<n$, $\textbf{X}_{(m)}\in \mathbb{R}^{m}$ and $\textbf{X}_{(n-m)}\in \mathbb{R}^{n-m}$, then the d.g. of $\textbf{X}_{(m)}$ is
$$\hat{g}_{m}(u)=\sum_{j=0}^{N-1}\frac{(N-1)_{j}}{j!} u^{N-1-j}\prod_{k=0}^{s_{1}}\int_{0}^{\infty}\frac{y^{\frac{n-m}{2}+j-1}e^{-\frac{a(s_{1})_{k}}{k!}u^{s_{1}-k}y^{k}}}
              {(1+e^{-b\sum_{l=0}^{s_{2}}\frac{(s_{2})_{l}}{l!}u^{s_{2}-l}y^{l}})^{2r}}dy ,$$
 where $(x)_{n}=x(x-1)(x-2)\cdots(x-n+1)$ and $\Phi_{2r}^{*}$ is the generalized Hurwitz-Lerch zeta function,
\end{flushleft}
\textbf{Proof.} By formula $(3.3)$, we have
\begin{equation}
\begin{split}
\hat{g}_{m}(u)&=\int_{u}^{\infty}(t-u)^{\frac{n-m}{2}-1}g_{n}(t)dt
         =\int_{0}^{\infty}y^{\frac{n-m}{2}-1}
         \frac{(y+u)^{N-1}e^{-a(y+u)^{s_{1}}}}{(1+e^{-b(y+u)^{s_{2}}})^{2r}}dy\nonumber\\
&=\int_{0}^{\infty}\sum_{j=0}^{N-1}\frac{(N-1)(N-2)\cdots(N-1-j+1)}{j!} u^{N-1-j}\frac{y^{\frac{n-m}{2}+j-1}e^{-a\sum_{k=0}^{s_{1}}\frac{(s_{1})_{k}}{k!}u^{s_{1}-k}y^{k}}}
              {(1+e^{-b\sum_{l=0}^{s_{2}}\frac{(s_{2})_{l}}{l!}u^{s_{2}-l}y^{l}})^{2r}}dy\\
&=\sum_{j=0}^{N-1}\frac{(N-1)_{j}}{j!} u^{N-1-j}\prod_{k=0}^{s_{1}}\int_{0}^{\infty}\frac{y^{\frac{n-m}{2}+j-1}e^{-\frac{a(s_{1})_{k}}{k!}u^{s_{1}-k}y^{k}}}
              {(1+e^{-b\sum_{l=0}^{s_{2}}\frac{(s_{2})_{l}}{l!}u^{s_{2}-l}y^{l}})^{2r}}dy.
\end{split}
\end{equation}
In addition, setting $s_{1}=s_{2}=1$,
\begin{equation}
\begin{split}
\hat{g}_{m}(u)
&=\sum_{j=0}^{N-1}\frac{(N-1)_{j}}{j!} e^{-au}u^{N-1-j}\int_{0}^{\infty}\frac{y^{\frac{n-m}{2}+j-1}e^{-ay}}{(1+e^{-bu}e^{-by})^{2r}}dy\nonumber\\
&=\sum_{j=0}^{N-1}\frac{(N-1)_{j}}{j!}u^{N-1-j}e^{-au}\frac{\Gamma(\frac{n-m}{2}+j)}{b^{\frac{n-m}{2}+j}}\Phi_{2r}^{*}(-e^{-bu},\frac{n-m}{2}+j,\frac{a}{b})\\
&=\sum_{j=0}^{N-1}\frac{(N-1)_{j}}{j!}\frac{\Gamma(\frac{n-m}{2}+j)}{b^{\frac{n-m}{2}+j}}u^{N-1-j}e^{-au}\Phi_{2r}^{*}(-e^{-bu},\frac{n-m}{2}+j,\frac{a}{b})
\end{split}
\end{equation}
follows.\\
Setting $s_{1}=s_{2}=s$, $r=0$,
 \begin{equation}
\begin{split}
\hat{g}_{m}(u)=&\sum_{j=0}^{N-1}\frac{(N-1)_{j}}{j!} u^{N-1-j}\int_{0}^{\infty}
y^{\frac{n-m}{2}+j-1}\exp\left(-\sum_{k=0}^{s}\frac{a(s)_{k}}{k!}y^{k}u^{s-k}\right)dy\nonumber\\
=&\sum_{j=0}^{N-1}\frac{(N-1)_{j}}{j!}u^{N-1-j}\int_{0}^{\infty}y^{\frac{n-m}{2}+j-1}
  \exp\left(-au^{s}-a\sum_{k=1}^{s}\frac{(s)_{k}}{k!}y^{k}u^{s-k}\right)dy\\
 &(letting ~x=u^{\frac{s-k}{k}}y)\\
  =&\sum_{j=0}^{N-1}\frac{(N-1)_{j}}{j!}u^{N-1-j}\exp(-au^{s})\prod_{k=1}^{s}\int_{0}^{\infty}
   \left(u^{-\frac{s-k}{k}}x\right)^{\frac{n-m}{2}+j-1}\exp\left(-a\frac{(s)_{k}}{k!}x^{k}\right)u^{-\frac{s-k}{k}}dy\\
  =&\sum_{j=0}^{N-1}\frac{(N-1)_{j}}{j!}u^{N-1-j}\exp(-au^{s})\prod_{k=1}^{s}\frac{1}{k}u^{-\frac{s-k}{k}(\frac{n-m}{2}+j)}\\
  &\times\int_{0}^{\infty}\left(\frac{k!}{a(s)_{k}}\right)^{\frac{1}{k}(\frac{n-m}{2}+j)-1}
  x^{\frac{1}{k}(\frac{n-m}{2}+j)}\exp(-x)\frac{k!}{r(s)_{k}}dx\\
  =&\sum_{j=0}^{N-1}\frac{(N-1)_{j}}{j!}u^{N-1-j}\exp(-au^{s})\prod_{k=1}^{s}\frac{1}{k}u^{-\frac{s-k}{k}(\frac{n-m}{2}+j)}
  \left(\frac{k!}{a(s)_{k}}\right)^{\frac{1}{k}(\frac{n-m}{2}+j)}\Gamma(\frac{1}{k}(\frac{n-m}{2}+j))\\
  =&\sum_{j=0}^{N-1}\prod_{k=1}^{s}\frac{(N-1)_{j}}{j!}\left(\frac{k!}{a(s)_{k}}\right)^{\frac{1}{k}(\frac{n-m}{2}+j)}
  \frac{\Gamma(\frac{1}{k}(\frac{n-m}{2}+j))}{k}u^{N-1-j-(1-\frac{s}{k})(\frac{n-m}{2}+j)}\exp(-au^{s})
\end{split}
\end{equation}
follows.
\par
The theorem above tells us that $g_{n}(\cdot),~\hat{g}_{m}(\cdot)$ are different in form rather than dimension changing.
\begin{flushleft}
\textbf{Corollary 3.1}
\end{flushleft}
\begin{enumerate}[1)]
  \item Let $\textbf{X}=(\textbf{X}_{(m)}^{\emph{T}},\textbf{X}_{(n-m)}^{\emph{T}})^{\emph{T}}\sim Lo_{n}(\bm{\mu},\mathbf{\Sigma},g_{n})$, where $g_{n}$ is defined as (2.5) then $$\widehat{g}_{m}(u)=\Gamma(\frac{n-m}{2})e^{-u}\Phi_{2}^{*}(-e^{-u},\frac{n-m}{2},1).$$
  \item Let $\textbf{X}=(\textbf{X}_{(m)}^{\emph{T}},\textbf{X}_{(n-m)}^{\emph{T}})^{\emph{T}}\sim Ko_{n}(\bm{\mu},\mathbf{\Sigma},g_{n})$, where $g_{n}$ is defined as (2.4) then
      $$\widehat{g}_{m}(u)=\sum_{j=0}^{N-1}\prod_{k=1}^{s_{1}}\frac{(N-1)_{j}}{j!}\left(\frac{k!}{a(s_{1})_{k}}\right)^{\frac{1}{k}(\frac{n-m}{2}+j)}
  \frac{\Gamma(\frac{1}{k}(\frac{n-m}{2}+j))}{k}u^{N-1-j-(1-\frac{s_{1}}{k})(\frac{n-m}{2}+j)}\exp(-au^{s_{1}}).$$
\end{enumerate}
\textbf{Proof.} We use Theorem 3.1 to conclude the density generator of $\textbf{X}_{(m)}$.
\begin{enumerate}[1)]
  \item Setting $N=a=b=r=1$, $s_{1}=s_{2}=1$ in Theorem 3.1 the result follows.
  \item Setting $s_{1}=s_{2}=s$, $r=0$, in Theorem 3.1 the result follows.
\end{enumerate}
\begin{flushleft}
\textbf{Remark 3.1} According to Fang et al. (1990), for $1\leq m\leq n-2$, the marginal density generators are related by $$g_{m+2}(x)=-\frac{1}{\pi}g_{m}'(x),~x>0~(a.e.),\eqno(3.4) $$ where $g_{m}'(\cdot)$ is the derivative of $g_{m}(\cdot)$. If $N\leq1$, (3.4) can be applied to the density generators of GL distribution.
\end{flushleft}
\section{Conditional distributions}
\par
Consider the partitions of $\textbf{X},~\bm{\mu},~\mathbf{\Sigma}$ as follows:
\par
$$\textbf{X}=\left(
\begin{matrix}
\textbf{X}^{(1)}\\ \textbf{X}^{(2)}
\end{matrix}
\right),~
\bm{\mu}=\left(
\begin{matrix}
\bm{\mu}^{(1)}\\ \bm{\mu}^{(2)}
\end{matrix}
\right),~
\mathbf{\Sigma}=\left(
\begin{matrix}
\mathbf{\Sigma}_{11}&\mathbf{\Sigma}_{12}\\ \mathbf{\Sigma}_{21}&\mathbf{\Sigma}_{22}
\end{matrix}
\right),
\eqno(4.1)$$
where $\textbf{X}^{(1)},~\bm{\mu}^{(1)}\in \mathbb{R}^{m}~(m<n),~\textbf{X}^{(2)},~\bm{\mu}^{(2)}\in \mathbb{R}^{n-m}$, $\mathbf{\Sigma}_{11}$ is an $m\times m$ matrix, $\mathbf{\Sigma}_{12}$ is an $m\times(n-m)$ matrix,
$\mathbf{\Sigma}_{21}$ is an $(n-m)\times m$ matrix and $\mathbf{\Sigma}_{22}$ is an $(n-m)\times(n-m)$ matrix.
\par
\begin{flushleft}
\textbf{Theorem 4.1} Let $\textbf{X}\sim GL_{n}(\bm{\mu},\mathbf{\Sigma},g_{n})$ where $g_{n}$ is defined as $(2.9)$. Conditionally on $\textbf{X}^{(2)}=\textbf{x}^{(2)}$, we have the conditional distribution of $\textbf{X}^{(1)}$.
\end{flushleft}
\begin{enumerate}[1)]
\item $\textbf{X}^{(1)}\sim Ell_{m}(\bm{\mu}_{1.2},\mathbf{\Sigma}_{11.2},g_{(1.2)})$ where $$\bm{\mu}_{1.2}=\bm{\mu}^{(1)}+\mathbf{\Sigma}_{12}\mathbf{\Sigma}_{22}^{-1}(\textbf{x}^{(2)}-\bm{\mu}^{(2)}),
      ~\mathbf{\Sigma}_{11.2}=\mathbf{\Sigma}_{11}-\mathbf{\Sigma}_{12}\mathbf{\Sigma}_{22}^{-1}\mathbf{\Sigma}_{21}.$$
\item The density generator of $\textbf{X}^{(1)}|\textbf{X}^{(2)}=\mathbf{x}^{(2)}$ can be written as
$$
g_{(1.2)}(t)
=\frac{(t+q(\textbf{x}^{(2)}))^{N-1}e^{-a\left(t+q(\textbf{x}^{(2)})\right)^{s_{1}}}}
{\left(1+e^{-b\left(t+q(\textbf{x}^{(2)})\right)^{s_{2}}}\right)^{2r}\sum_{j=0}^{N-1}\frac{(N-1)_{j}}{j!} q(\textbf{x}^{(2)})^{N-1-j}\prod_{k=0}^{s_{1}}I_{k,j}},
$$
where $q(\textbf{x}^{(2)})=(\textbf{x}^{(2)}-\bm{\mu}^{(2)})^{'}\mathbf{\Sigma}_{22}^{-1}(\textbf{x}^{(2)}-\bm{\mu}^{(2)})$,
$$I_{k,j}\triangleq I_{k,j}(m,a,b,s_{1},s_{2},r)=\int_{0}^{\infty}\frac{y^{\frac{m}{2}+j-1}e^{-\frac{a(s_{1})_{k}}{k!}q(\textbf{x}^{(2)})^{s_{1}-k}y^{k}}}
              {(1+e^{-b\sum_{l=0}^{s_{2}}\frac{(s_{2})_{l}}{l!}q(\textbf{x}^{(2)})^{s_{2}-l}y^{l}})^{2r}}dy.$$
\end{enumerate}
\par
\textbf{Proof.}
\begin{enumerate}[1)]
\item The result follows by Fang et al. (1990).
\item  Define an $n$-dimension vector $$\textbf{Y}=\left(
\begin{matrix}
\textbf{Y}^{(1)}\\ \textbf{Y}^{(2)}
\end{matrix}
\right)=\textbf{D}\textbf{X},~\textbf{D}=
\left(
\begin{matrix}
\textbf{I}&\-\mathbf{\Sigma}_{12}\mathbf{\Sigma}_{22}^{-1}\\0&\textbf{I}
\end{matrix}
\right),$$ where $\textbf{I}$ is an $m\times m$ identity matrix, $\mathbf{\Sigma}_{12},\mathbf{\Sigma}_{22}$ are given in $(4.1)$, to make \\$\textbf{Y}^{(1)}\in \mathbb{R}^{m}$, $\textbf{Y}^{(2)}\in \mathbb{R}^{n-m}$ are independent with each other. We have $$\textbf{Y}^{(1)}=\textbf{X}^{(1)}-\mathbf{\Sigma}_{12}\mathbf{\Sigma}_{22}^{-1}\textbf{X}^{(2)},
  \textbf{Y}^{(2)}=\textbf{X}^{(2)};$$
\begin{equation}
\begin{split}
\bm{\mu}_{\textbf{Y}}=&\textbf{D}\bm{\mu}=\left(
\begin{matrix}
\textbf{I}&\-\mathbf{\Sigma}_{12}\mathbf{\Sigma}_{22}^{-1}\\0&\textbf{I}
\end{matrix}
\right)
\left(
\begin{matrix}
\bm{\mu}^{(1)}\\ \bm{\mu}^{(2)}
\end{matrix}
\right)=
\left(
\begin{matrix}
\bar{\bm{\mu}}_{(1.2)}\\ \bm{\mu}^{(2)}
\end{matrix}
\right),
~\bar{\bm{\mu}}_{(1.2)}=\bm{\mu}^{(1)}-\mathbf{\Sigma}_{12}\mathbf{\Sigma}_{22}^{-1}\bm{\mu}^{(2)};\\
\mathbf{\Sigma}_{\textbf{Y}}=&\textbf{D}\mathbf{\Sigma}\textbf{D}^{'}\nonumber\\
=&\left(
\begin{matrix}
\textbf{I}&\-\mathbf{\Sigma}_{12}\mathbf{\Sigma}_{22}^{-1}\\0&\textbf{I}
\end{matrix}
\right)
\left(
\begin{matrix}
\mathbf{\Sigma}_{11}&\mathbf{\Sigma}_{12}\\ \mathbf{\Sigma}_{21}&\mathbf{\Sigma}_{22}
\end{matrix}
\right)
\left(
\begin{matrix}
\textbf{I}&0\\-\mathbf{\Sigma}_{12}\mathbf{\Sigma}_{22}^{-1}&\textbf{I}
\end{matrix}
\right)\\=&
\left(
\begin{matrix}
\mathbf{\Sigma}_{11}-\mathbf{\Sigma}_{12}\mathbf{\Sigma}_{22}^{-1}\mathbf{\Sigma}_{21}&0\\ \mathbf{\Sigma}_{21}&\mathbf{\Sigma}_{22}
\end{matrix}
\right)
\left(
\begin{matrix}
\textbf{I}&0\\-\mathbf{\Sigma}_{12}\mathbf{\Sigma}_{22}^{-1}&\textbf{I}
\end{matrix}
\right)\\
=&\left(
\begin{matrix}
\mathbf{\Sigma}_{11.2}&0\\0&\mathbf{\Sigma}_{22}
\end{matrix}
\right),
\end{split}
\end{equation}
where $\mathbf{\Sigma}_{11.2}=\mathbf{\Sigma}_{11}-\mathbf{\Sigma}_{12}\mathbf{\Sigma}_{22}^{-1}\mathbf{\Sigma}_{21}$.
Then $\textbf{Y}\sim GL_{n}(\bm{\mu}_{\textbf{Y}},\mathbf{\Sigma}_{\textbf{Y}},g_{n})$, the d.g. of \textbf{Y} is same as \textbf{X}'s d.g.. Moreover, the density generators of $\textbf{X}^{(1)}$, $\textbf{X}^{(2)}$ are identical with density generators of $\textbf{Y}^{(1)}$, $\textbf{Y}^{(2)}$ respectively, i.e. $\hat{g}_{\textbf{X}^{(1)}}(t)=\hat{g}_{\textbf{Y}^{(1)}}(t)=\hat{g}_{m}(t)$, $\tilde{g}_{\textbf{X}^{(2)}}(u)=\tilde{g}_{\textbf{Y}^{(2)}}(u)=\tilde{g}_{n-m}(u).$ $g_{n}(\cdot),~\hat{g}_{m}(\cdot),~\tilde{g}_{n-m}(\cdot)$ are pairwise different in form rather than dimension changing.
\begin{equation}
\begin{split}
f_{\textbf{Y}}=&
\frac{C_{n}}{\sqrt{|\mathbf{\Sigma}_{\textbf{Y}}|}}
g_{n}\left((\textbf{y}-\bm{\mu}_{\textbf{Y}})^{\emph{T}}\mathbf{\Sigma}_{\textbf{Y}}^{-1}(\textbf{y}-\bm{\mu}_{\textbf{Y}})\right)\nonumber\\
=&\frac{C_{n}}{\sqrt{|\mathbf{\Sigma}_{\textbf{Y}}|}}
g_{n}\left((\textbf{Y}^{(1)}-\bar{\bm{\mu}}_{1.2})^{\emph{T}}\mathbf{\Sigma}_{11.2}^{-1}(\textbf{Y}^{(1)}-\bar{\bm{\mu}}_{1.2})
+(\textbf{Y}^{(2)}-\bm{\mu}^{(2)})^{\emph{T}}\mathbf{\Sigma}_{22}^{-1}(\textbf{Y}^{(2)}-\bm{\mu}^{(2)})\right)\\
=&\frac{C_{n}}{\sqrt{|\mathbf{\Sigma}_{\textbf{Y}}|}}
g_{n}\left((\textbf{X}^{(1)}-\mathbf{\Sigma}_{12}\mathbf{\Sigma}_{22}^{-1}\textbf{X}^{(2)}-\bar{\bm{\mu}}_{1.2})^{\emph{T}}
      \mathbf{\Sigma}_{11.2}^{-1}(\textbf{X}^{(1)}-\mathbf{\Sigma}_{12}\mathbf{\Sigma}_{22}^{-1}\textbf{X}^{(2)}-\bar{\bm{\mu}}_{1.2})\right.\\
& \left.+(\textbf{X}^{(2)}-\bm{\mu}^{(2)})^{\emph{T}}\mathbf{\Sigma}_{22}^{-1}(\textbf{X}^{(2)}-\bm{\mu}^{(2)})\right)\\
=&\frac{C_{n}}{\sqrt{|\mathbf{\Sigma}_{\textbf{Y}}|}}
   g_{n}\left((\textbf{X}^{(1)}-\bm{\mu}_{1.2})^{\emph{T}}\mathbf{\Sigma}_{11.2}^{-1}(\textbf{X}^{(1)}-\bm{\mu}_{1.2})
+(\textbf{X}^{(2)}-\bm{\mu}^{(2)})^{\emph{T}}\mathbf{\Sigma}_{22}^{-1}(\textbf{X}^{(2)}-\bm{\mu}^{(2)})\right),\\
\bm{\mu}_{1.2}=&\bm{\mu}^{(1)}+\mathbf{\Sigma}_{12}\mathbf{\Sigma}_{22}^{-1}(\textbf{X}^{(2)}-\bm{\mu}^{(2)}),\\
f_{\textbf{Y}^{(2)}}
=&\frac{\tilde{C}_{n-m}}
       {\sqrt{|\mathbf{\Sigma}_{22}|}}
  \tilde{g}_{n-m}\left((\textbf{X}^{(2)}-\bm{\mu}^{(2)})^{\emph{T}}\mathbf{\Sigma}_{22}^{-1}(\textbf{X}^{(2)}-\bm{\mu}^{(2)})\right),\\
f_{\textbf{Y}^{(1)}|\textbf{Y}^{(2)}}
=&\frac{\hat{C}_{m}}{\sqrt{|\mathbf{\Sigma}_{11.2}|}}
  g_{(1.2)}\left((\textbf{X}^{(1)}-\bm{\mu}_{1.2})^{\emph{T}}\mathbf{\Sigma}_{11.2}^{-1}(\textbf{X}^{(1)}-\bm{\mu}_{1.2})\right).
\end{split}
\end{equation}
Since $$ f_{\textbf{Y}^{(1)}|\textbf{Y}^{(2)}}=\frac{f_{\textbf{Y}}}{f_{\textbf{Y}^{(2)}}},$$
we have $$g_{(1.2)}\left((\textbf{X}^{(1)}-\bm{\mu}_{1.2})^{\emph{T}}\mathbf{\Sigma}_{11.2}^{-1}(\textbf{X}^{(1)}-\bm{\mu}_{1.2})\right)
  =\frac{g_{n}\left((\textbf{X}^{(1)}-\bm{\mu}_{1.2})^{\emph{T}}\mathbf{\Sigma}_{11.2}^{-1}(\textbf{X}^{(1)}-\bm{\mu}_{1.2})
                  +q(\textbf{x}^{(2)})\right)}{\tilde{g}_{(n-m)}(q(\textbf{x}^{(2)}))},$$
then we have
$$
g_{(1.2)}(t)
=\frac{(t+q(\textbf{x}^{(2)}))^{N-1}e^{-a\left(t+q(\textbf{x}^{(2)})\right)^{s_{1}}}}
{\left(1+e^{-b\left(t+q(\textbf{x}^{(2)})\right)^{s_{2}}}\right)^{2r}\sum_{j=0}^{N-1}\frac{(N-1)_{j}}{j!} q(\textbf{x}^{(2)})^{N-1-j}\prod_{k=0}^{s_{1}}I_{k,j}},
$$
and $$\bm{\mu}_{1.2}=\bm{\mu}^{(1)}+\mathbf{\Sigma}_{12}\mathbf{\Sigma}_{22}^{-1}(\textbf{x}^{(2)}-\bm{\mu}^{(2)}),
      ~\mathbf{\Sigma}_{11.2}=\mathbf{\Sigma}_{11}-\mathbf{\Sigma}_{12}\mathbf{\Sigma}_{22}^{-1}\mathbf{\Sigma}_{21}.$$
where $q(\textbf{x}^{(2)})=(\textbf{x}^{(2)}-\bm{\mu}^{(2)})^{\emph{T}}\mathbf{\Sigma}_{22}^{-1}(\textbf{x}^{(2)}-\bm{\mu}^{(2)})$,
$$I_{k,j}\triangleq I_{k,j}(m,a,b,s_{1},s_{2},r)=\int_{0}^{\infty}\frac{y^{\frac{m}{2}+j-1}e^{-\frac{a(s_{1})_{k}}{k!}q(\textbf{x}^{(2)})^{s_{1}-k}y^{k}}}
              {(1+e^{-b\sum_{l=0}^{s_{2}}\frac{(s_{2})_{l}}{l!}q(\textbf{x}^{(2)})^{s_{2}-l}y^{l}})^{2r}}dy.$$
In addition, setting $s_{1}=s_{2}=1$,
$$
g_{(1.2)}(t)=\frac{(t+q(\textbf{x}^{(2)}))^{N-1}e^{-at}}
                   {(1+e^{-b(t+q(\textbf{x}^{(2)}))})^{2r}\sum_{j=0}^{N-1}\frac{(N-1)_{j}}{j!}\frac{\Gamma(\frac{m}{2})+j}{b^{\frac{m}{2}+j}}(q(\textbf{x}^{(2)}))^{N-1-j}
             \Phi_{2r}^{*}(-e^{-bq(\textbf{x}^{(2)})},\frac{m}{2}+j,\frac{a}{b})}
$$ holds,
where $q(\textbf{x}^{(2)})=(\textbf{x}^{(2)}-\bm{\mu}^{(2)})^{\emph{T}}\mathbf{\Sigma}_{22}^{-1}(\textbf{x}^{(2)}-\bm{\mu}^{(2)})$.
\end{enumerate}

\begin{flushleft}
\textbf{Corollary 4.1}
\end{flushleft}
\begin{enumerate}[1)]
  \item Supposing $\textbf{X}\sim Lo_{n}(\bm{\mu},\mathbf{\Sigma},g_{n})$ where $g_{n}$ is written as (2.5), the partitions of $\textbf{X},~\bm{\mu},~\mathbf{\Sigma}$ are same as $(4.1)$. The density generator of $\textbf{X}^{(1)}$ conditionally on $\textbf{X}^{(2)}=\textbf{x}^{(2)}$ is
      $$g_{(1.2)}(t)=\left[\Gamma(\frac{m}{2})\Phi_{2}^{*}(-e^{-q(\textbf{x}^{(2)})},\frac{m}{2},1)\right]^{-1}
      \frac{e^{-t}}{(1+e^{-(t+q(\textbf{x}^{(2)}))})^{2}}.$$
  \item Supposing $\textbf{X}\sim Ko_{n}(\bm{\mu},\mathbf{\Sigma},g_{n})$ where $g_{n}$ is written as (2.4), the partitions of $\textbf{X},~ \bm{\mu},~ \mathbf{\Sigma}$ is same as $(4.1)$. The density generator of $\textbf{X}^{(1)}$ conditionally on $\textbf{X}^{(2)}=\textbf{x}^{(2)}$ is
      $$g_{(1.2)}(t)=\frac{\left(t+q(\textbf{x}^{(2)})\right)^{N-1}e^{-a(t+q(\textbf{x}^{(2)}))^{s_{1}}}}{\sum_{j=0}^{N-1}\prod_{k=1}^{s_{1}}
      \gamma_{k,j}t^{N-1-j-(1-\frac{s_{1}}{k})(\frac{n-m}{2}+j)}e^{-at^{s_{1}}}},$$
      where $$\gamma_{x,y}\triangleq\gamma_{x,y}(N,a,s_{1},n,m,b)=\frac{(N-1)_{y}}{y!}\left(\frac{x!}{a(s_{1})_{x}}\right)^{\frac{1}{x}(\frac{n-m}{2}+y)}
  \frac{\Gamma(\frac{1}{x}(\frac{n-m}{2}+y))}{x}.$$
  \end{enumerate}
\section{Characteristic functions and characteristic generators}
\begin{flushleft}
\textbf{Theorem 5.1}  Let $\textbf{X}\sim GL_{n}(\bm{\mu},\mathbf{\Sigma},g_{n})$ where $g_{n}$ is defined as $(2.9)$. The characteristic function of $\textbf{X}$ can be expressed as follows.
\end{flushleft}
\begin{enumerate}[1)]
  \item If $n=1$,
$$\psi_{\textbf{X}}(t)=C_{1}e^{it\mu}\int_{0}^{\infty}\frac{x^{N-\frac{3}{2}}e^{-ax^{s}}}{(1+e^{-bx^{s}})^{2r}}\cos(t\sigma\sqrt{x})dx, ~t\in(-\infty, \infty).$$
  \item If $n>1$,
$$
\psi_{\textbf{X}}(\textbf{t})=e^{i\textbf{t}^{\emph{T}}\bm{\mu}}\sum_{j=0}^{\infty}\frac{(-1)^{j}}{(2j)!}\frac{q_{j}}{q_{0}}
\frac{\Phi_{2r}^{*}(-1,\frac{1}{s}(\frac{n}{2}+j+N-1),\frac{a}{b})}{\Phi_{2r}^{*}(-1,\frac{1}{s}(\frac{n}{2}+N-1),\frac{a}{b})}
(\textbf{t}^{\emph{T}}\mathbf{\Sigma}\textbf{t})^{j},
$$
where $\Phi_{2r}^{*}$ is the generalized Hurwitz-Lerch zeta function, $B(\cdot)$ is the Beta function, $\textbf{t}=(t_{1}, t_{2},\cdots, t_{n})^{\emph{T}},~t_{i}\in(-\infty, \infty),$
$$q_{x}\triangleq q_{x}(N,n,b,s)=\frac{\Gamma(\frac{1}{s}(\frac{n}{2}+x+N-1))}{b^{\frac{1}{s}(\frac{n}{2}+x+N-1)}}B(\frac{n-1}{2},\frac{2x+1}{2}).$$
\end{enumerate}
\textbf{Proof.}
If $n=1$,
\begin{equation}
\begin{split}
\psi_{\textbf{X}}(t)=&E(e^{itx})=\int_{-\infty}^{\infty}e^{itx}\frac{C_{1}}{\sigma}g_{1}[(\frac{x-\mu}{\sigma})^{2}]dx\nonumber\\
=&\int_{-\infty}^{\infty}C_{1}e^{it\mu}e^{it\sigma x}g_{1}(x^{2})dx
=C_{1}e^{it\mu}\int_{-\infty}^{\infty}[\cos(t\sigma x)+i\sin(t\sigma x)]g_{1}(x^{2})dx\\
=&C_{1}e^{it\mu}\int_{0}^{\infty}\cos(t\sigma \sqrt{y})y^{-\frac{1}{2}}g_{1}(y)dy
=C_{1}e^{it\mu}\int_{0}^{\infty}\cos(t\sigma \sqrt{y})\frac{y^{N-\frac{3}{2}}e^{-ay^{s}}}{(1+e^{-by^{s}})^{2r}}dy.\\
\end{split}
\end{equation}
If $n>1$,
\begin{equation}
\begin{split}
\psi_{\textbf{X}}(\textbf{t})=& E(e^{i\textbf{t}^{\emph{T}}\textbf{X}})=e^{i\textbf{t}^{\emph{T}}\bm{\mu}}E(e^{i\textbf{t}^{\emph{T}}R\mathbf{\Sigma}^{\frac{1}{2}}\textbf{U}^{(n)}})
                       =e^{i\textbf{t}^{\emph{T}}\bm{\mu}}E[\Omega_{n}(R^{2}\textbf{t}^{\emph{T}}\mathbf{\Sigma}\textbf{t})]\nonumber\\
                       =&e^{i\textbf{t}^{\emph{T}}\bm{\mu}}\int_{0}^{\infty}\Omega_{n}(v\textbf{t}^{\emph{T}}\mathbf{\Sigma}\textbf{t})
                        \frac{1}{\int_{0}^{\infty}t^{\frac{n}{2}-1}g_{n}(t)dt}v^{\frac{n}{2}-1}\frac{v^{N-1}\exp(-av^{s})}{(1+\exp(-bv^{s}))^{2r}}dv\\
                       =&e^{i\textbf{t}^{\emph{T}}\bm{\mu}}I_{1},
\end{split}
\end{equation}
where $\Omega_{n}(\|t\|^{2}), t\in \mathbb{R}^{n}$ is the characteristic function of $\mathbf{U^{(n)}}$ (Fang et al. (1990))
$$\Omega_{n}(\|t\|^{2})=\frac{1}{B(\frac{n-1}{2},\frac{1}{2})}\int_{0}^{\pi}\exp(i\|t\|^{2}\cos\theta)\sin^{n-2}\theta d\theta,$$
\begin{equation}
\begin{split}
I_{1}=&\int_{0}^{\infty}\Omega_{n}(v\textbf{t}^{\emph{T}}\mathbf{\Sigma}\textbf{t})\frac{1}{\int_{0}^{\infty}t^{\frac{n}{2}-1}g_{n}(t)dt}
        \frac{v^{\frac{n}{2}+N-2}\exp(-av^{s})}{(1+\exp(-bv^{s}))^{2r}}dv\nonumber\\
=&\int_{0}^{\infty}\frac{1}{B(\frac{n-1}{2},\frac{1}{2})}
  \left[\frac{\Gamma(\frac{1}{s}(\frac{n}{2}+N-1))}{b^{\frac{1}{s}(\frac{n}{2}+N-1)}s}\Phi_{2r}^{*}(-1,\frac{1}{s}(\frac{n}{2}+N-1),\frac{a}{b})\right]^{-1}\\
 &\times\int_{0}^{\pi}\exp(iv^{\frac{1}{2}}(\textbf{t}^{\emph{T}}\mathbf{\Sigma}\textbf{t})^{\frac{1}{2}}\cos\theta)\sin^{n-2}\theta
   d\theta\frac{v^{\frac{n}{2}+N-2}\exp(-av^{s})}{(1+\exp(-bv^{s}))^{2r}}dv\\
=&\frac{1}{B(\frac{n-1}{2},\frac{1}{2})}
  \left[\frac{\Gamma(\frac{1}{s}(\frac{n}{2}+N-1))}{b^{\frac{1}{s}(\frac{n}{2}+N-1)}s}\Phi_{2r}^{*}(-1,\frac{1}{s}(\frac{n}{2}+N-1),\frac{a}{b})\right]^{-1}
   \int_{0}^{\pi}I_{2}\sin^{n-2}\theta d\theta,\\
I_{2}=&\int_{0}^{\infty}\frac{v^{\frac{n}{2}+N-2}\exp(-av^{s}+iv^{\frac{1}{2}}(\textbf{t}^{\emph{T}}\mathbf{\Sigma}\textbf{t})^{\frac{1}{2}}\cos\theta)}
                       {(1+\exp(-bv^{s}))^{2r}}dv\\ \nonumber
      =&\sum_{l=0}^{\infty}\frac{i^{l}(\textbf{t}^{\emph{T}}\mathbf{\Sigma}\textbf{t})^{\frac{l}{2}}\cos^{l}\theta}{l!}
  \int_{0}^{\infty}\frac{v^{\frac{n+l}{2}+N-2}\exp(-av^{s})}{(1+\exp(-bv^{s}))^{2r}}dv\\
   =&\sum_{l=0}^{\infty}\frac{i^{l}(\textbf{t}^{\emph{T}}\mathbf{\Sigma}\textbf{t})^{\frac{l}{2}}\cos^{l}\theta}{l!}
  \int_{0}^{\infty}\frac{1}{s}\frac{v^{\frac{1}{s}(\frac{n+l}{2}+N-1)-1}\exp(-av)}{(1+\exp(-bv))^{2r}}dv\\
=&\sum_{l=0}^{\infty}\frac{\Gamma(\frac{1}{s}(\frac{n+l}{2}+N-1))}{b^{\frac{1}{s}(\frac{n+l}{2}+N-1)}s}
   \frac{i^{l}(\textbf{t}^{\emph{T}}\mathbf{\Sigma}\textbf{t})^{\frac{l}{2}}\cos^{l}\theta}{l!}\Phi_{2r}^{*}(-1,\frac{1}{s}(\frac{n+l}{2}+N-1),\frac{a}{b}).
\end{split}
\end{equation}
Then, we have
\begin{equation}
\begin{split}
I_{1}=&\frac{1}{B(\frac{n-1}{2},\frac{1}{2})}\left[\frac{\Gamma(\frac{1}{s}(\frac{n}{2}+N-1))}{b^{\frac{1}{s}(\frac{n}{2}+N-1)}s}\Phi_{2r}^{*}(-1,\frac{1}{s}(\frac{n}{2}+N-1),\frac{a}{b})\right]^{-1}\\
   &\times\sum_{l=0}^{\infty}\frac{\Gamma(\frac{1}{s}(\frac{n+l}{2}+N-1))}{b^{\frac{1}{s}(\frac{n+l}{2}+N-1)}s}
   \frac{i^{l}(t^{'}\mathbf{\Sigma}t)^{\frac{l}{2}}}{l!}\Phi_{2r}^{*}(-1,\frac{1}{s}(\frac{n+l}{2}+N-1),\frac{a}{b})
   \int_{0}^{\pi}\sin^{n-2}\theta\cos^{l}\theta d\theta\\
=&\frac{1}{B(\frac{n-1}{2},\frac{1}{2})}\left[\frac{\Gamma(\frac{1}{s}(\frac{n}{2}+N-1))}{b^{\frac{1}{s}(\frac{n}{2}+N-1)}}\Phi_{2r}^{*}(-1,\frac{1}{s}(\frac{n}{2}+N-1),\frac{a}{b})\right]^{-1}\nonumber\\
 &\times[\sum_{j=0}^{\infty}\frac{\Gamma(\frac{1}{s}(\frac{n+2j+1}{2}+N-1))}{b^{\frac{1}{s}(\frac{n+2j+1}{2}+N-1)}}
   \frac{i^{2j+1}(\textbf{t}^{\emph{T}}\mathbf{\Sigma}\textbf{t})^{\frac{2j+1}{2}}}{(2j+1)!}\Phi_{2r}^{*}(-1,\frac{1}{s}(\frac{n+2j+1}{2}+N-1),\frac{a}{b})\cdot0\\
 &+\sum_{j=0}^{\infty}\frac{\Gamma(\frac{1}{s}(\frac{n+2j}{2}+N-1))}{b^{\frac{1}{s}(\frac{n+2j}{2}+N-1)}}
   \frac{i^{2j}(\textbf{t}^{\emph{T}}\mathbf{\Sigma}\textbf{t})^{\frac{2j}{2}}}{(2j)!}\Phi_{2r}^{*}(-1,\frac{1}{s}(\frac{n+2j}{2}+N-1),\frac{a}{b})B(\frac{n-1}{2},\frac{2j+1}{2})]\\
=&\left[B(\frac{n-1}{2},\frac{1}{2})\frac{\Gamma(\frac{1}{s}(\frac{n}{2}+N-1))}{b^{\frac{1}{s}(\frac{n}{2}+N-1)}}\Phi_{2r}^{*}(-1,\frac{1}{s}(\frac{n}{2}+N-1),\frac{a}{b})\right]^{-1}\\
 &\times\sum_{j=0}^{\infty}B(\frac{n-1}{2},\frac{2j+1}{2})\frac{\Gamma(\frac{1}{s}(\frac{n+2j}{2}+N-1))}{b^{\frac{1}{s}(\frac{n+2j}{2}+N-1)}}
   \frac{(-1)^{j}(\textbf{t}^{\emph{T}}\mathbf{\Sigma}\textbf{t})^{j}}{(2j)!}\Phi_{2r}^{*}(-1,\frac{1}{s}(\frac{n+2j}{2}+N-1),\frac{a}{b}).
\end{split}
\end{equation}
Finally, we obtain
$$
\psi_{\textbf{X}}(\textbf{t})=e^{it^{\emph{T}}\mathbf{\mu}}\sum_{j=0}^{\infty}\frac{(-1)^{j}}{(2j)!}\frac{q_{j}}{q_{0}}
\frac{\Phi_{2r}^{*}(-1,\frac{1}{s}(\frac{n}{2}+j+N-1),\frac{a}{b})}{\Phi_{2r}^{*}(-1,\frac{1}{s}(\frac{n}{2}+N-1),\frac{a}{b})}(\textbf{t}^{\emph{T}}\mathbf{\Sigma}\textbf{t})^{j},$$
where $$q_{x}\triangleq q_{x}(N,n,b,s)=\frac{\Gamma(\frac{1}{s}(\frac{n}{2}+x+N-1))}{b^{\frac{1}{s}(\frac{n}{2}+x+N-1)}}B(\frac{n-1}{2},\frac{2x+1}{2}).$$
\begin{flushleft}
\textbf{Corollary 5.1} Let $\textbf{X}\sim GL_{n}(\bm{\mu},\mathbf{\Sigma},g_{n})$ where $g_{n}$ is defined as $(2.9)$. The characteristic generator of $\textbf{X}$ can be expressed as follows.
\end{flushleft}
\begin{enumerate}[1)]
  \item  If $n=1$,
  $$\phi_{\textbf{X}}(t)=C_{1}\int_{0}^{\infty}\frac{x^{N-\frac{3}{2}}e^{-ax^{s}}}{(1+e^{-bx^{s}})^{2r}}\cos(t\sigma\sqrt{x})dx,
  ~t\in(-\infty, \infty).$$
  \item If $n>1$, letting $\textbf{u}_{n}=(u_{1},u_{2},\cdots,u_{n})^{\emph{T}}$,
  $$\phi_{\textbf{X}}(\|\textbf{u}_{n}\|^{2})
 =\sum_{j=0}^{\infty}\frac{(-1)^{j}}{(2j)!}\frac{q_{j}}{q_{0}}
\frac{\Phi_{2r}^{*}(-1,\frac{1}{s}(\frac{n}{2}+j+N-1),\frac{a}{b})}{\Phi_{2r}^{*}(-1,\frac{1}{s}(\frac{n}{2}+N-1),\frac{a}{b})}\|\textbf{u}_{n}\|^{2j},
$$
where $$q_{x}\triangleq q_{x}(N,n,b,s)=\frac{\Gamma(\frac{1}{s}(\frac{n}{2}+x+N-1))}{b^{\frac{1}{s}(\frac{n}{2}+x+N-1)}}B(\frac{n-1}{2},\frac{2x+1}{2}).$$
\end{enumerate}
\begin{flushleft}
\textbf{Proof.} The results directly follow by the definition of characteristic generator.
\end{flushleft}
\begin{flushleft}
\textbf{Remark 5.1} $\Omega_{n}(\|t\|^{2})$ can be expressed in the following alternative forms:
\end{flushleft}
$$\Omega_{n}(\|t\|^{2})=\frac{\Gamma(\frac{n}{2})}{\sqrt{\pi}}\sum_{k=0}^{\infty}\frac{(-1)^{k}\|t\|^{2k}}{(2k)!}
   \frac{\Gamma(\frac{2k+1}{2})}{\Gamma(\frac{n+2k}{2})},$$
$$\Omega_{n}(\|t\|^{2})=_{0}F_{1}(\frac{n}{2};-\frac{1}{4}\|t\|^{2}).$$
We can obtain the following equivalent forms of the characteristic functions and characteristic generators with dimension $n>1$.
$$\psi_{\textbf{X}}(t)=e^{i\textbf{t}^{\emph{T}}\bm{\mu}}\sum_{k=0}^{\infty}\frac{(-1)^{k}}{(2k)!}\gamma_{k}(N,n,b,s)
 \frac{\Phi_{2r}^{*}(-1,\frac{1}{s}(\frac{n}{2}+N+k-1),\frac{a}{b})}{\Phi_{2r}^{*}(-1,\frac{1}{s}(\frac{n}{2}+N-1),\frac{a}{b})}
 (\textbf{t}^{\emph{T}}\mathbf{\Sigma}\textbf{t})^{k}, \eqno{(5.1)}$$
$$\phi_{\textbf{X}}(\|\textbf{u}_{n}\|^{2})
 =\sum_{k=0}^{\infty}\frac{(-1)^{k}}{(2k)!}\gamma_{k}(N,n,b,s)
 \frac{\Phi_{2r}^{*}(-1,\frac{1}{s}(\frac{n}{2}+N+k-1),\frac{a}{b})}{\Phi_{2r}^{*}(-1,\frac{1}{s}(\frac{n}{2}+N-1),\frac{a}{b})}
 (\|\textbf{u}_{n}\|)^{2k},\eqno{(5.2)}$$
where $$\gamma_{k}(N,n,b,s)=\frac{\Gamma(\frac{n}{2})}{\pi^{\frac{1}{2}}b^{\frac{1}{s}(\frac{n}{2}+N-1)}}
\frac{\Gamma(k+\frac{1}{2})\Gamma(\frac{1}{s}(\frac{n}{2}+N+k-1))}{\Gamma(k+\frac{n}{2})\Gamma(\frac{1}{s}(\frac{n}{2}+N-1))}.$$
$$
\psi_{\textbf{X}}(\textbf{t})=e^{i\textbf{t}^{\emph{T}}\bm{\mu}}\sum_{k=0}^{\infty}
\frac{\Gamma(\frac{1}{s}(\frac{n}{2}+N+k-1))}{\Gamma(\frac{1}{s}(\frac{n}{2}+N-1))b^{\frac{k}{s}}4^{k}(\frac{n}{2})^{[k]}k!}
\frac{\Phi_{2r}^{*}(-1,\frac{1}{s}(\frac{n}{2}+N+k-1),\frac{a}{b})}{\Phi_{2r}^{*}(-1,\frac{1}{s}(\frac{n}{2}+N-1),\frac{a}{b})}
(\textbf{t}^{\emph{T}}\Sigma \textbf{t})^{k},\eqno{(5.3)}$$
$$\phi_{\textbf{X}}(\|\textbf{u}_{n}\|^{2})
 =\sum_{k=0}^{\infty}
\frac{\Gamma(\frac{1}{s}(\frac{n}{2}+N+k-1))}{\Gamma(\frac{1}{s}(\frac{n}{2}+N-1))b^{\frac{k}{s}}4^{k}(\frac{n}{2})^{[k]}k!}
\frac{\Phi_{2r}^{*}(-1,\frac{1}{s}(\frac{n}{2}+N+k-1),\frac{a}{b})}{\Phi_{2r}^{*}(-1,\frac{1}{s}(\frac{n}{2}+N-1),\frac{a}{b})}
(\|\textbf{u}_{n}\|)^{2k}.\eqno{(5.4)}$$
\par
  Similar as the density generator, the characteristic generator (c.g.) of the class of multivariate elliptically symmetric
distribution defined in Definition 2.1 is not dimensionally coherent. In other word, the characteristic function $\psi_{\textbf{X}}(\cdot)$ and the characteristic generator $\phi_{\textbf{X}}(\cdot)$ are related to the dimension of \textbf{X}. More details will be discussed in Section 7.
\section{Moments}
\begin{flushleft}
\textbf{Theorem 6.1} Let \textbf{X} $\sim GBL_{n}(\bm{\mu},\mathbf{\Sigma}, g_{n})$ where $g_{n}$ is defined as $(2.9)$.
\begin{enumerate}[1)]
\item The expectation and the covariance are:
      $$E(\textbf{X})=\bm{\mu},
      ~Cov(\textbf{X})=\frac{1}{n}\frac{\Gamma(\frac{1}{s}(N+\frac{n}{2}))\Phi_{2r}^{*}(-1,\frac{1}{s}(N+\frac{n}{2}),\frac{a}{b})}
       {b^{\frac{1}{s}}\Gamma(\frac{1}{s}(N+\frac{n}{2}-1))\Phi_{2r}^{*}(-1,\frac{1}{s}(N+\frac{n}{2}-1),\frac{a}{b})}\mathbf{\Sigma};$$
\item For any integers $m_{1},\cdots,m_{n}$, with $m=\sum_{i=1}^{n}m_{i},$ the product moments of $\textbf{Z}:=\mathbf{\Sigma}^{-\frac{1}{2}}(\textbf{X}-\bm{\mu})$ are
    $$E(\prod_{i=1}^{n}Z_{i}^{m_{i}})=\frac{1}{n}\frac{\Gamma(\frac{1}{s}(N+\frac{n}{2}+\frac{m}{2}-1))\Phi_{2r}^{*}(-1,\frac{1}{s}(N+\frac{n}{2}+\frac{m}{2}-1),\frac{a}{b})}
       {(\frac{n}{2})^{[l]}b^{\frac{m}{2s}}\Gamma(\frac{1}{s}(N+\frac{n}{2}-1))\Phi_{2r}^{*}(-1,\frac{1}{s}(N+\frac{n}{2}-1),\frac{a}{b})}
       \prod_{i=1}^{n}\frac{(2l_{i})!}{4^{l_{i}}(l_{i})!},$$
    where $x^{[n]}=x(x+1)\cdots(x+n-1)$ and $\Phi_{2r}^{*}$ is the generalized Hurwitz-Lerch zeta function,
    if $m_{i}=2l_{i}$ are even, $i=1,\cdots,n, m=2l$;
    $$E(\prod_{i=1}^{n}Z_{i}^{m_{i}})=0,$$ if at least one of the $m_{i}$ is odd.
\end{enumerate}
\end{flushleft}
\par
\textbf{Proof.} According to $(1.2)$ we have for real number $p>0$,
\begin{equation}
\begin{split}
E(R^{p})=&\frac{1}{\int_{0}^{\infty}t^{n-1}g_{n}(t^{2})dt}\int_{0}^{\infty}z^{n-1+p}g_{n}(z^{2})dz\\
        =&\left[\int_{0}^{\infty}t^{n-1}\frac{t^{2(N-1)}e^{-at^{2s}}}{(1+e^{-bt^{2s}})^{2r}}dt\right]^{-1}
          \int_{0}^{\infty}\frac{z^{2N+n+p-3}e^{-az^{2s}}}{(1+e^{-bz^{2s}})^{2r}}dz\nonumber\\
&(setting~x=bz^{2s})\\
=&\frac{2b^{\frac{1}{s}(N+\frac{n}{2}-1)}s}{\Gamma(\frac{1}{s}(N+\frac{n}{2}-1))}
        \left[\Phi_{2r}^{*}(-1,\frac{1}{s}(N+\frac{n}{2}-1),\frac{a}{b})\right]^{-1}\\
&\times\int_{0}^{\infty}\frac{1}{2s}b^{-\frac{1}{2s}(2N+n+p-2)}\frac{x^{\frac{1}{2s}(2N+n+p-2)-1}e^{-\frac{a}{b}x}}{(1+e^{-x})^{2r}}dx\\
=&\frac{\Gamma(\frac{1}{s}(N+\frac{n}{2}+\frac{p}{2}-1))\Phi_{2r}^{*}(-1,\frac{1}{s}(N+\frac{n}{2}+\frac{p}{2}-1),\frac{a}{b})}
       {b^{\frac{p}{2s}}\Gamma(\frac{1}{s}(N+\frac{n}{2}-1))\Phi_{2r}^{*}(-1,\frac{1}{s}(N+\frac{n}{2}-1),\frac{a}{b})}.
\end{split}
\end{equation}
\begin{enumerate}[1)]
\item Since $\textbf{X}=\bm{\mu}+R\textbf{A}^{\emph{T}}\textbf{U}^{(n)}$ and $E(\textbf{U}^{(n)})=0$, we have $$E(\textbf{X})=\bm{\mu}+E(R)\textbf{A}^{\emph{T}}E(\textbf{U}^{(n)})=\bm{\mu},$$ and
\begin{equation}
\begin{split}
Cov(\textbf{X})&=Cov(R\textbf{A}^{\emph{T}}\textbf{U}^{(n)})
                   =E(R^{2})\textbf{A}^{\emph{T}}Cov(\textbf{U}^{(n)})\textbf{A}
                   =\frac{1}{n}E(R^{2})\mathbf{\Sigma}\nonumber\\
&=\frac{1}{n}\frac{\Gamma(\frac{1}{s}(N+\frac{n}{2}))\Phi_{2r}^{*}(-1,\frac{1}{s}(N+\frac{n}{2}),\frac{a}{b})}
       {b^{\frac{1}{s}}\Gamma(\frac{1}{s}(N+\frac{n}{2}-1))\Phi_{2r}^{*}(-1,\frac{1}{s}(N+\frac{n}{2}-1),\frac{a}{b})}\mathbf{\Sigma}.
\end{split}
\end{equation}
\item By Eqs.(2.18) and (3.6) in Fang et al. (1990),
$$ E(\prod_{i=1}^{n} Z_{i}^{m_{i}})
                                   =E(R^{m})E(\prod_{i=1}^{n}u_{i}^{m_{i}}),$$
where $E(\prod_{i=1}^{n}u_{i}^{m_{i}})=\frac{1}{(\frac{n}{2})^{[l]}}
                                        \prod_{i=1}^{n}\frac{(2l_{i})!}{4^{l_{i}}(l_{i})!}$,
    if $m_{i}=2l_{i}~(i=1,\cdots,n)$ are even, $m=2l$; $E(\prod_{i=1}^{n}u_{i}^{m_{i}})=0$, if at least one of the $m_{i}$ is odd.
\end{enumerate}
\begin{flushleft}
\textbf{Corollary 6.1}
\begin{enumerate}[1)]
  \item Let \textbf{X} $\sim Ko_{n}(\mathbf{\mu},\mathbf{\Sigma}, g_{n})$ where $g_{n}$ is defined as $(2.4)$.
\begin{enumerate}[(1)]
\item The expectation and the covariance are:
 $$E(\textbf{X})=\bm{\mu},~Cov(\textbf{X})=\frac{1}{n}\frac{\Gamma(\frac{1}{s_{1}}(N+\frac{n}{2}))}{a^{\frac{1}{s_{1}}}\Gamma(\frac{1}{s_{1}}(N+\frac{n}{2}-1))}
      \mathbf{\Sigma};$$
\item For any integers $m_{1},\cdots,m_{n}$, with $m=\sum_{i=1}^{n}m_{i},$ the product moments of $\textbf{Z}:=\mathbf{\Sigma}^{-\frac{1}{2}}(\textbf{X}-\bm{\mu})$ are
    $$E(\prod_{i=1}^{n}Z_{i}^{m_{i}})=\frac{\Gamma(\frac{1}{s_{1}}(N+\frac{n}{2}+\frac{m}{2}-1))}
    {(\frac{n}{2})^{[l]}a^{\frac{m}{2s_{1}}}\Gamma(\frac{1}{s_{1}}(N+\frac{n}{2}-1))}\prod_{i=1}^{n}\frac{(2l_{i})!}{4^{l_{i}}(l_{i})!},$$
    where $x^{[n]}=x(x+1)\cdots(x+n-1)$, $\Phi_{2r}^{*}$ is the generalized Hurwitz-Lerch zeta function,
     if $m_{i}=2l_{i}~(i=1,\cdots,n)$ are even, $m=2l$;
    $$E(\prod_{i=1}^{n}Z_{i}^{m_{i}})=0,$$ if at least one of the $m_{i}$ is odd.
    \end{enumerate}
  \item Let \textbf{X} $\sim Lo_{n}(\mathbf{\mu},\mathbf{\Sigma}, g_{n})$ where $g_{n}$ is defined as $(2.5)$.
  \begin{enumerate}[(1)]
  \item The expectation and the covariance are:
      $$E(\textbf{X})=\mathbf{\mu},~Cov(\textbf{X})=\frac{1}{n}\frac{\Gamma(1+\frac{n}{2})\Phi_{2}^{*}(-1,N+\frac{n}{2},1)}
       {\Gamma(\frac{n}{2})\Phi_{2}^{*}(-1,\frac{n}{2},1)}\mathbf{\Sigma};$$
\item For any integers $m_{1},\cdots,m_{n}$, with $m=\sum_{i=1}^{n}m_{i},$ the product moments of $Z:=\mathbf{\Sigma}^{-\frac{1}{2}}(\textbf{X}-\mathbf{\mu})$ are
    $$E(\prod_{i=1}^{n}Z_{i}^{m_{i}})=\frac{\Gamma(\frac{n}{2}+\frac{m}{2})\Phi_{2}^{*}(-1,\frac{n}{2}+\frac{m}{2},1)}{(\frac{n}{2})^{[l]}\Gamma(\frac{n}{2})\Phi_{2}^{*}(-1,\frac{n}{2},1)}\prod_{i=1}^{n}\frac{(2l_{i})!}{4^{l_{i}}(l_{i})!},$$
     if $m_{i}=2l_{i}~(i=1,\cdots,n)$ are even, $m=2l$; $E(\prod_{i=1}^{n}u_{i}^{m_{i}})=0$, if at least one of the $m_{i}$ is odd,
    where $x^{[n]}=x(x+1)\cdots(x+n-1)$, $\Phi_{2r}^{*}$ is the generalized Hurwitz-Lerch zeta function.
\end{enumerate}
\end{enumerate}
\end{flushleft}
\section{Linear transform and marginal distributions}
\par
It has been mentioned in lots of researches that if $\textbf{X}\sim Ell_{n}(\bm{\mu},\mathbf{\Sigma},\phi)$, rank($\mathbf{\Sigma}$)$=k$, \\ \textbf{B} is an $n\times m$ matrix and \textbf{v} is an $m\times 1$ vector, then
$$\textbf{v}+\textbf{B}^{\emph{T}}\textbf{X}\sim Ell_{m}(\textbf{v}+\textbf{B}^{\emph{T}}\bm{\mu},\textbf{B}^{\emph{T}}\mathbf{\Sigma}\textbf{B},\phi).$$In the theorem above the characteristic generators of elliptical distributions are regarded as unrelated to dimension i.e. the characteristic generator would not change with the dimension during the liner transform. However the dimension coherent property is inapplicable for the characteristic generator of GL distribution. We will demonstrate GL distribution's liner transform property in the following theorem.
\par
\begin{flushleft}
\textbf{Theorem 7.1} Assuming $\textbf{X}\sim GL_{n}(\bm{\mu},\mathbf{\Sigma},\phi)$ with stochastic representation $\textbf{X}=\bm{\mu}+R\textbf{A}^{\emph{T}}\textbf{U}^{(n)}$, $\textbf{Y}=\textbf{B}\textbf{X}+\textbf{b}$, where \textbf{B} is an $m\times n$ matrix, $1\leq m<n$, rank$(\textbf{B})=m$ and $\textbf{b}\in \mathbb{R}^{n}$, the d.g. of \textbf{X} is
$$g_{n}(x)=\frac{x^{N-1}\exp(-ax)}{(1+\exp(-bx))^{2r}}.$$
\end{flushleft}
\begin{enumerate}[1)]
\item $\textbf{Y}\sim Ell_{m}(\textbf{B}\bm{\mu}+\textbf{b},\textbf{B}\mathbf{\Sigma}\textbf{B}^{\emph{T}},g_{(m),y})$, where
      $$g_{(m),y}=\sum_{j=0}^{N-1}\frac{(N-1)_{j}}{j!}
        \frac{\Gamma(\frac{n-m}{2}+j)}{b^{\frac{n-m}{2}+j}}u^{N-1-j}e^{-au}\Phi_{2r}^{*}(-e^{-bu},\frac{n-m}{2},\frac{a}{b}).$$
\item $\textbf{Y}\sim Ell_{m}(\textbf{B}\bm{\mu}+\textbf{b},\textbf{B}\mathbf{\Sigma}\textbf{B}^{\emph{T}},\phi_{(m),y})$, where
    $$\phi_{(m),y}(t)=\sum_{j=0}^{N-1}\frac{(N-1)_{j}}{j!}
                            \frac{\Gamma(\frac{n-1}{2}+j)}{b^{\frac{n-1}{2}+j}}
                            \Phi_{2r}^{*}(-e^{-bu},\frac{n-1}{2}+j,\frac{a}{b}) \int_{0}^{\infty}y^{N-j-\frac{3}{2}}e^{-ay}\cos(t\sigma_{\textbf{Y}}\sqrt{y})dy,$$
$t\in(-\infty, \infty)$, if $m=1$;
$$\phi_{(m),y}(\|\bm{\xi}_{(m)}\|^{2})
 =\frac{\sum_{j=0}^{N-1}\sum_{l=0}^{\infty}\sum_{\omega=0}^{\infty}\frac{(-1)^{\omega}}{(2\omega)!}\alpha_{l}^{*}\beta_{j}
 q_{\omega}^{*}(N,m,a,b,j,l)}{\sum_{k=0}^{N-1}\sum_{j=0}^{\infty}\alpha_{j}^{*}\beta_{k}q_{0}^{*}(N,m,a,1,j,l)}\|\bm{\xi}_{(m)}\|^{2\omega},
 ~\xi_{i}\in(-\infty, \infty),$$
$$\alpha_{x}^{*}\triangleq\alpha_{x}^{*}(N,n,m,a,b,r)=\frac{(-1)^{x}\Gamma(2r+x)}{x!(x+\frac{a}{b})^{\frac{n}{2}+N-1}},~
\beta_{x}\triangleq\beta_{x}(N,n,m,b,r)=\frac{(N-1)_{x}}{x!}\frac{\Gamma(\frac{n-m}{2}+x)}{b^{\frac{n-m}{2}+x}\Gamma(2r)},$$
$$q_{\omega}^{*}(N,m,a,b,j,l)=
\frac{\Gamma(\frac{m}{2}+\omega+N-j-1)B(\frac{m-1}{2},\frac{2\omega+1}{2})}{b^{\frac{m}{2}+N+\omega-j-1}(l+\frac{a}{b})^{\omega}},~ \omega=0, 1, 2, \cdots,$$
$\bm{\xi}_{(m)}=(\mathbf{\xi}_{1},\mathbf{\xi}_{2},\cdots,\mathbf{\xi}_{m})'$,
if $m>1$.
\end{enumerate}
\par
\begin{flushleft}
\textbf{Proof.} $\textbf{Y}$ is no longer GL distributed but elliptically distributed still.
\end{flushleft}
\begin{enumerate}[1)]
\item Applying Theorem 3.1 the result follows.
\item When $m$=1,
\begin{equation}
\begin{split}
\psi_{\textbf{Y}}(t)=&E(e^{it\textbf{Y}})
=\int_{-\infty}^{\infty}e^{ity}
\frac{\tilde{C}_{1}}{\sigma_{\textbf{Y}}}\tilde{g}_{1}[(\frac{y-\mu_{\textbf{Y}}}{\sigma_{\textbf{Y}}})^{2}]dx\nonumber
=\tilde{C}_{1}e^{it\mu_{\textbf{Y}}}\int_{0}^{\infty}\cos(t\sigma_{\textbf{Y}}\sqrt{y})\tilde{g}_{1}(y)y^{-\frac{1}{2}}dy\\
=&\tilde{C}_{1}e^{it\mu_{\textbf{Y}}}\int_{0}^{\infty}\cos(t\sigma_{\textbf{Y}}\sqrt{y})
  \sum_{j=0}^{N-1}\frac{(N-1)_{j}}{j!}\frac{\Gamma(\frac{n-1}{2}+j)}{b^{\frac{n-1}{2}+j}}
  y^{N-j-\frac{3}{2}}e^{-ay}\Phi_{2r}^{*}(-e^{-by},\frac{n-1}{2}+j,\frac{a}{b})dy\\
=&\tilde{C}_{1}e^{it\mu_{\textbf{Y}}}\sum_{j=0}^{N-1}\frac{(N-1)_{j}}{j!}\frac{\Gamma(\frac{n-1}{2}+j)}{b^{\frac{n-1}{2}+j}}
  \int_{0}^{\infty}y^{N-j-\frac{3}{2}}e^{-ay}\cos(t\sigma_{\textbf{Y}}\sqrt{y})
  \Phi_{2r}^{*}(-e^{-by},\frac{n-1}{2}+j,\frac{a}{b})dy,
  \end{split}
\end{equation}
therefore,
$$\phi_{\textbf{Y}}(t^{2})=\sum_{j=0}^{N-1}C_{N-1}^{j}
                            \frac{\Gamma(\frac{n-1}{2}+j)}{b^{\frac{n-1}{2}+j}}
                            \Phi_{2r}^{*}(-e^{-bu},\frac{n-1}{2}+j,\frac{a}{b}) \int_{0}^{\infty}y^{N-j-\frac{3}{2}}e^{-ay}\cos(t\sigma_{\textbf{Y}}\sqrt{y})dy.$$

When $m>1$,
\begin{equation}
\begin{split}
\psi_{\textbf{Y}}(\textbf{t})=&E(e^{i\textbf{t}^{\emph{T}}\textbf{Y}})
=e^{i\textbf{t}^{\emph{T}}\bm{\mu}_{\textbf{Y}}}\frac{1}{\int_{0}^{\infty}t^{\frac{n}{2}-1}g_{n}(t)dt}\nonumber\\
&\times\int_{0}^{\infty}\Omega_{n}(v\textbf{t}^{\emph{T}}\mathbf{\Sigma}_{\textbf{Y}}\textbf{t})v^{\frac{m}{2}-1}
\sum_{j=0}^{N-1}\frac{(N-1)_{j}}{j!}\frac{\Gamma(\frac{n-m}{2}+j)}{b^{\frac{n-m}{2}+j}}v^{N-1-j}e^{-av}
\Phi_{2r}^{*}(-e^{-bv},\frac{n-m}{2},\frac{a}{b})dv\\
=&e^{i\textbf{t}^{\emph{T}}\bm{\mu}_{\textbf{Y}}}I_{1}^{-1}\frac{1}{B(\frac{m-1}{2},\frac{1}{2})}\int_{0}^{\pi}\sin^{m-2}\theta I_{2}d\theta,
\end{split}
\end{equation}
where
\begin{equation}
\begin{split}
I_{1}^{-1}=&\frac{1}{\int_{0}^{\infty}t^{\frac{m}{2}-1}g_{m}(t)dt}\nonumber\\
=&\left[\int_{0}^{\infty}t^{\frac{m}{2}-1}\sum_{k=0}^{N-1}\frac{(N-1)_{k}}{k!}\frac{\Gamma(\frac{n-m}{2}+k)}{b^{\frac{n-m}{2}+k}}
  t^{N-1-k}e^{-at}\Phi_{2r}^{*}(-e^{-bt},\frac{n-m}{2}+k,\frac{a}{b})dt\right]^{-1}\\
=&\left[\sum_{k=0}^{N-1}\frac{(N-1)_{k}}{k!}\frac{\Gamma(\frac{n-m}{2}+k)}{b^{\frac{n-m}{2}+k}}\frac{1}{\Gamma(2r)}
  \int_{0}^{\infty}t^{N+\frac{m}{2}-k-2}e^{-at}\sum_{j=0}^{\infty}
  \frac{\Gamma(2r+j)}{j!}\frac{(-1)^{j}e^{-bjt}}{(j+\frac{a}{b})^{\frac{n-m}{2}+k}}dt\right]^{-1}\\
=&\left[\sum_{k=0}^{N-1}\beta_{k}\sum_{j=0}^{\infty}\alpha_{j}\frac{1}{(j+\frac{a}{b})^{\frac{n-m}{2}+k}}
  \int_{0}^{\infty}t^{N+\frac{m}{2}-k-2}e^{-(a+bj)t}dt\right]^{-1}\\
=&\left[\sum_{k=0}^{N-1}\sum_{j=0}^{\infty}\beta_{k}\alpha_{j}\frac{1}{(j+\frac{a}{b})^{\frac{n-m}{2}+k}}
\frac{1}{(a+bj)^{\frac{m}{2}+N-1-k}}\int_{0}^{\infty}t^{N+\frac{m}{2}-k-2}e^{-t}dt\right]^{-1}\\
=&\left[\sum_{k=0}^{N-1}\sum_{j=0}^{\infty}\alpha_{j}\beta_{k}\frac{\Gamma(\frac{m}{2}+N-k-1)}{(j+\frac{a}{b})^{\frac{n}{2}+N-1}}\right]^{-1}\\
=&\left[\sum_{k=0}^{N-1}\sum_{j=0}^{\infty}\alpha_{j}^{*}\beta_{k}\Gamma(\frac{m}{2}+N-k-1)\right]^{-1},
\end{split}
\end{equation}
$$\beta_{x}\triangleq\beta_{x}(N,n,m,b,r)=\frac{(N-1)_{x}}{x!}\frac{\Gamma(\frac{n-m}{2}+x)}{b^{\frac{n-m}{2}+x}\Gamma(2r)},$$
$$\alpha_{x}(r)\triangleq\alpha_{x}(r)=\frac{\Gamma(2r+x)(-1)^{x}}{x!},~\alpha_{x}^{*}\triangleq\alpha_{x}^{*}(N,n,a,b,r)=\frac{\alpha_{x}}{(x+\frac{a}{b})^{\frac{n}{2}+N-1}},
$$
\begin{equation}
\begin{split}
I_{2}=&\int_{0}^{\infty}\exp(iv^{\frac{1}{2}}(\textbf{t}^{\emph{T}}\mathbf{\Sigma_{\textbf{Y}}}\textbf{t})^{\frac{1}{2}}\cos\theta)\\
&\times\sum_{j=0}^{N-1}\frac{(N-1)_{j}}{j!}\frac{\Gamma(\frac{n-m}{2}+j)}{b^{\frac{n-m}{2}+j}}v^{\frac{m}{2}+N-2-j}e^{-av}
\Phi_{2r}^{*}(-e^{-bv},\frac{n-m}{2}+j,\frac{a}{b})dv\nonumber\\
=&\sum_{j=0}^{N-1}\frac{(N-1)_{j}}{j!}\frac{\Gamma(\frac{n-m}{2}+j)}{b^{\frac{n-m}{2}+j}}
\int_{0}^{\infty}v^{N+\frac{m}{2}-j-2}e^{iv^{\frac{1}{2}}(\textbf{t}^{\emph{T}}\mathbf{\Sigma}\textbf{t})^{\frac{1}{2}}\cos\theta-av}
\Phi_{2r}^{*}(-e^{-bv},\frac{n-m}{2}+j,\frac{a}{b})dv\\
=&\sum_{j=0}^{N-1}\beta_{j}
  \int_{0}^{\infty}v^{N+\frac{m}{2}-j-2}e^{iv^{\frac{1}{2}}(\textbf{t}^{\emph{T}}\mathbf{\Sigma_{\textbf{Y}}}\textbf{t})^{\frac{1}{2}}\cos\theta-av}
  \sum_{l=0}^{\infty}\frac{\Gamma(2r+l)}{l!}\frac{(-1)^{l}e^{-blv}}{(l+\frac{a}{b})^{\frac{n-m}{2}+j}}dv\\
=&\sum_{j=0}^{N-1}\sum_{l=0}^{\infty}\beta_{j}\alpha_{l}\frac{1}{(l+\frac{a}{b})^{\frac{n-m}{2}+j}}
\int_{0}^{\infty}v^{N+\frac{m}{2}-j-2}e^{-av-blv}\sum_{q=0}^{\infty}\frac{i^{q}v^{\frac{q}{2}}}{q!}
 (\textbf{t}^{\emph{T}}\mathbf{\Sigma_{\textbf{Y}}}\textbf{t})^{\frac{q}{2}}\cos^{q}\theta dv\\
 =&\sum_{j=0}^{N-1}\sum_{l=0}^{\infty}\sum_{q=0}^{\infty}\beta_{j}\alpha_{l}
 \frac{(\textbf{t}^{\emph{T}}\mathbf{\Sigma_{\textbf{Y}}}\textbf{t})^{\frac{q}{2}}\cos^{q}\theta}{(l+\frac{a}{b})^{\frac{n-m}{2}+j}}\frac{i^{q}}{q!}
 \int_{0}^{\infty}v^{N+\frac{q}{2}+\frac{m}{2}-j-2}e^{-av-blv}dv\\
=&\sum_{j=0}^{N-1}\sum_{l=0}^{\infty}\sum_{q=0}^{\infty}\beta_{j}\alpha_{l}
 \frac{(\textbf{t}^{\emph{T}}\mathbf{\Sigma_{\textbf{Y}}}\textbf{t})^{\frac{q}{2}}\cos^{q}\theta}{(l+\frac{a}{b})^{\frac{n-m}{2}+j}}\frac{i^{q}}{q!}
 \frac{\Gamma(\frac{m+q}{2}+N-k-1)}{(a+bl)^{\frac{m}{2}+\frac{q}{2}+N-j-1}}.
 \end{split}
\end{equation}
 Then we have
 \begin{equation}
\begin{split}
&\int_{0}^{\pi}\sin^{m-2}\theta I_{2}d\theta\\
=&\sum_{j=0}^{N-1}\sum_{l=0}^{\infty}\sum_{q=0}^{\infty}\alpha_{l}\beta_{j}
 \frac{(\textbf{t}^{\emph{T}}\mathbf{\Sigma_{\textbf{Y}}}\textbf{t})^{\frac{q}{2}}}{(l+\frac{a}{b})^{\frac{n-m}{2}+j}}\frac{i^{q}}{q!}
 \frac{\Gamma(\frac{m+q}{2}+N-j-1)}{(a+bl)^{\frac{m}{2}+\frac{q}{2}+N-j-1}}\int_{0}^{\pi}\cos^{q}\theta\sin^{m-2}\theta d\theta\nonumber\\
 =&\sum_{j=0}^{N-1}\sum_{l=0}^{\infty}\sum_{\omega=0}^{\infty}\frac{(-1)^{\omega}}{(2\omega)!}\alpha_{l}^{*}\beta_{j}
 \frac{\Gamma(\frac{m}{2}+\omega+N-j-1)B(\frac{m-1}{2},\frac{2\omega+1}{2})}{b^{\frac{m}{2}+N+\omega-j-1}(l+\frac{a}{b})^{\omega}}
 (\textbf{t}^{\emph{T}}\mathbf{\Sigma_{\textbf{Y}}}\textbf{t})^{\omega},
 \end{split}
\end{equation}
 \begin{equation}
\begin{split}
\psi_{\textbf{Y}}(\textbf{t})=&e^{i\textbf{t}^{\emph{T}}\bm{\mu}_{\textbf{Y}}}
\frac{\sum_{j=0}^{N-1}\sum_{l=0}^{\infty}\sum_{\omega=0}^{\infty}\frac{(-1)^{\omega}}{(2\omega)!}\alpha_{l}^{*}\beta_{j}
 \frac{\Gamma(\frac{m}{2}+\omega+N-j-1)B(\frac{m-1}{2},\frac{2\omega+1}{2})}{b^{\frac{m}{2}+N+\omega-j-1}(l+\frac{a}{b})^{\omega}}}
 {\sum_{k=0}^{N-1}\sum_{j=0}^{\infty}\alpha_{j}^{*}\beta_{k}\Gamma(\frac{m}{2}+N-k-1)B(\frac{m-1}{2},\frac{1}{2})}
 (\textbf{t}^{\emph{T}}\mathbf{\Sigma_{\textbf{Y}}}\textbf{t})^{\omega}\nonumber\\
 =&e^{i\textbf{t}^{\emph{T}}\bm{\mu}_{\textbf{Y}}}
 \frac{\sum_{j=0}^{N-1}\sum_{l=0}^{\infty}\sum_{\omega=0}^{\infty}\frac{(-1)^{\omega}}{(2\omega)!}\alpha_{l}^{*}\beta_{j}
 q_{\omega}^{*}(N,m,a,b,j,l)}
 {\sum_{k=0}^{N-1}\sum_{j=0}^{\infty}\alpha_{j}^{*}\beta_{k}q_{0}^{*}(N,m,a,1,j,l)}
 (\textbf{t}^{\emph{T}}\mathbf{\Sigma_{\textbf{Y}}}\textbf{t})^{\omega},
  \end{split}
\end{equation}
where $$ q_{\omega}^{*}(N,m,a,b,j,l)=
\frac{\Gamma(\frac{m}{2}+\omega+N-j-1)B(\frac{m-1}{2},\frac{2\omega+1}{2})}{b^{\frac{m}{2}+N+\omega-j-1}(l+\frac{a}{b})^{\omega}}, ~\omega=0, 1, 2, \cdots.$$
 Therefore,
 $$\phi_{(m),y}(\|\bm{\xi}_{(m)}\|^{2})
 =\frac{\sum_{j=0}^{N-1}\sum_{l=0}^{\infty}\sum_{\omega=0}^{\infty}\frac{(-1)^{\omega}}{(2\omega)!}\alpha_{l}^{*}\beta_{j}
 q_{\omega}^{*}(N,m,a,b,j,l)}
 {\sum_{k=0}^{N-1}\sum_{j=0}^{\infty}\alpha_{j}^{*}\beta_{k}q_{0}^{*}(N,m,a,b,j,l)}\|\bm{\xi}_{(m)}\|^{2\omega}.
 $$
 \end{enumerate}
\begin{flushleft}
\textbf{Remark 7.1} Similar as Remark 5.1, we can obtain the following equivalent forms of the characteristic functions and characteristic generators of $\textbf{Y}=\textbf{B}\textbf{X}+\textbf{b}$ with dimension $m>1$.
\end{flushleft}
$$\psi_{\textbf{Y}}(\textbf{t})
=e^{i\textbf{t}^{\emph{T}}\bm{\mu}_{\textbf{Y}}}\frac{\sum_{l=0}^{\infty}\sum_{k=0}^{N-1} \sum_{p=0}^{\infty}
\frac{(-1)^{l}}{(2l)!}\alpha_{p}^{*}\beta_{k}A_{l}\Gamma(\frac{m}{2}+N-k+l-1)}
{\sum_{k=0}^{\infty}\sum_{j=0}^{N-1}\alpha_{k}^{*}\beta_{j}\Gamma(\frac{m}{2}+N-j-1)}(\textbf{t}^{\emph{T}}\mathbf{\Sigma_{\textbf{Y}}}\textbf{t})^{l},
\eqno{(7.1)}$$
 $$\phi_{(m),y}(\|\bm{\xi}_{(m)}\|^{2})
=\frac{\sum_{l=0}^{\infty}\sum_{k=0}^{N-1}\sum_{p=0}^{\infty}\frac{(-1)^{l}}{(2l)!}\alpha_{p}^{*}\beta_{k}A_{l}\Gamma(\frac{m}{2}+N-k+l-1)}
{\sum_{k=0}^{\infty}\sum_{j=0}^{N-1}\alpha_{k}^{*}\beta_{j}\Gamma(\frac{m}{2}+N-j-1)}\|\bm{\xi}_{(m)}\|^{2l};
\eqno{(7.2)}$$
$$\psi_{\textbf{Y}}(\textbf{t})=e^{i\textbf{t}^{\emph{T}}\bm{\mu}_{\textbf{Y}}}\frac{\sum_{l=0}^{\infty}\sum_{j=0}^{N-1}\sum_{k=0}^{\infty}\alpha_{k}\beta_{j}^{*}
\Phi_{2r}^{*}(-e^{-bv},\frac{n-m}{2}+j,\frac{a}{b})}{\sum_{k=0}^{\infty}\sum_{j=0}^{N-1}\alpha_{k}^{*}\beta_{j}\Gamma(\frac{m}{2}+N-j-1)}
\frac{\Gamma(\frac{m}{2})}{\pi^{\frac{1}{2}}4^{l}(\frac{m}{2})^{[l]}}(\textbf{t}^{\emph{T}}\mathbf{\Sigma}_{\textbf{Y}}\textbf{t})^{l},
\eqno{(7.3)}$$
$$\phi_{(m),y}(\|\bm{\xi}_{(m)}\|^{2})=\frac{\sum_{l=0}^{\infty}\sum_{j=0}^{N-1}\sum_{k=0}^{\infty}\alpha_{k}\beta_{j}^{*}
\Phi_{2r}^{*}(-e^{-bv},\frac{n-m}{2}+j,\frac{a}{b})}{\sum_{k=0}^{\infty}\sum_{j=0}^{N-1}\alpha_{k}^{*}\beta_{j}\Gamma(\frac{m}{2}+N-j-1)}
\frac{\Gamma(\frac{m}{2})}{\pi^{\frac{1}{2}}4^{l}(\frac{m}{2})^{[l]}}\|\bm{\xi}_{(m)}\|^{2l},\eqno{(7.4)}$$
 where
 $$A_{x}\triangleq A_{x}(N,m,a,b,p,k)=\frac{B(\frac{m}{2},x+\frac{1}{2})}{B(\frac{m}{2}+x,\frac{1}{2})b^{N-k+x}(p+\frac{a}{b})^{x-\frac{m}{2}+1}},$$
  $$\beta_{x}\triangleq\beta_{x}(N,n,m,r)=\frac{(N-1)_{x}}{x!}\frac{\Gamma(\frac{n-m}{2}+x)}{\Gamma(2r)b^{\frac{n-m}{2}+x}},~
  \beta_{x}^{*}\triangleq\beta_{x}^{*}(N,n,m,b,r)=\frac{\beta_{x}\Gamma(N+\frac{m}{2}-x-1)}{b^{N+\frac{m}{2}-x-1}},$$
  $$\alpha_{x}\triangleq\alpha_{x}(r)=\frac{(-1)^{x}\Gamma(2r+x)}{x!},~
  \alpha_{x}^{*}\triangleq\alpha_{x}^{*}(N,m,a,b,r)=\frac{\alpha_{x}}{(x+\frac{a}{b})^{\frac{n}{2}+N-1}}.$$
Here $\Phi_{2r}^{*}$ is the generalized Hurwitz-Lerch zeta function.
\begin{flushleft}
\textbf{Corollary 7.1} (\emph{\textbf{Marginal distributions}}) Supposing $\textbf{X}\sim GL_{n}(\bm{\mu},\mathbf{\Sigma},g_{n})$, where $g_{n}$ is defined as $(2.9)$ with $s_{1}=s_{2}=1$, the partitions of $\textbf{X},~\bm{\mu},~\mathbf{\Sigma}$ are given in (4.1), then
$$\textbf{X}^{(1)}\sim Ell_{m}(\bm{\mu}^{(1)},\mathbf{\Sigma}_{11},\widehat{g}_{m}),~
\textbf{X}^{(2)}\sim Ell_{n-m}(\bm{\mu}^{(2)},\mathbf{\Sigma}_{22},\widehat{g}_{n-m}),$$
where $$\widehat{g}_{m}(u)=\sum_{j=0}^{N-1}\frac{(N-1)_{j}}{j!}\frac{\Gamma(\frac{n-m}{2}+j)}{b^{\frac{n-m}{2}+j}}u^{N-1-j}e^{-au}\Phi_{2r}^{*}(-e^{-bu},\frac{n-m}{2}+j,\frac{a}{b}),$$
$\Phi_{2r}^{*}$ is the generalized Hurwitz-Lerch zeta function.
\end{flushleft}
\par
\textbf{Proof.} Taking $\textbf{B}_{1}=(\textbf{I}_{m},\textbf{0}_{m\times(n-m)}),~ \textbf{B}_{2}=(\textbf{0}_{m\times(n-m)},\textbf{I}_{n-m})$, $\textbf{X}^{(1)}=\textbf{B}_{1}\textbf{X}$,\\  $\textbf{X}^{(2)}=\textbf{B}_{2}\textbf{X}$ in Theorem 4.1 the result follows.
\par
On the basis of Theorem 5.1 with $s=1$ and Theorem 7.1, it is clear that the c.g. of GL distributed random vector depends on dimension.
\par
\begin{flushleft}
\textbf{Example} (\textbf{\emph{Local dependence function}}) Bairamov et al (2003) presented the local dependence function denoted by $H(x,y)$ based on regression concepts as follows:
\end{flushleft}
$$H(x,y)=\frac{E\{(X-E(X|Y=y))(Y-E(Y|X=x))\}}
{\sqrt{E\{(X-E(X|Y=y))^{2}\}}\sqrt{E\{(Y-E(Y|X=x))^{2}\}}}.$$
 Alternative representations of $H(x,y)$ are
$$H(x,y)=\frac{Cov(X,Y)+\xi_{Y}(x)\xi_{X}(y)}{\sqrt{Var(X)+\xi^{2}_{X}(y)}\sqrt{Var(Y)+\xi^{2}_{Y}(x)}},$$
$$H(x,y)=\frac{\rho+\phi_{X}(y)\phi_{Y}(x)}{\sqrt{1+\phi_{X}^{2}(y)}\sqrt{1+\phi_{Y}^{2}(x)}},$$
where$$\xi_{X}(y)=E(X|Y=y)-E(X),~\xi_{Y}(x)=E(Y|X=x)-E(Y),~\rho=\frac{Cov(X,Y)}{\sqrt{Var(X)}\sqrt{Var(Y)}},$$
$$\phi_{X}(y)=\frac{\xi_{X}(y)}{\sqrt{Var(X)}},~\phi_{Y}(x)=\frac{\xi_{Y}(x)}{\sqrt{Var(Y)}}.$$
This function can characterize the dependence structure of two random variables X, Y localized at the fixed point.
Suppose $\textbf{W}=(X,Y)^{\emph{T}}\sim GL_{2}(\bm{\mu},\mathbf{\Sigma},g_{2})$, where $g_{2}$ is defined as $(2.9)$ with $s=1$, $\widehat{I}=\int_{0}^{\infty}t^{2}\widehat{g}_{1}(t^{2})dt$. Without loss of generality, let
$$\bm{\mu}=\left(
\begin{matrix}
0\\0
\end{matrix}
\right),~
\mathbf{\Sigma}=\left(
\begin{matrix}
1&\rho^{'}\\ \rho^{'}&1
\end{matrix}
\right).
$$
We have
$$Cov(\textbf{W})
=\frac{N}{2b}\frac{\Phi_{2r}^{*}(-1,N+1,\frac{a}{b})}{\Phi_{2r}^{*}(-1,N,\frac{a}{b})}\rho',$$
$$E(Y|X=x)=\xi_{Y}(x)=\rho'x,~E(X|Y=y)=\xi_{X}(y)=\rho'y,$$ respectively.
$$Var(X)=\frac{2C_{1}}{\sqrt{1-\rho'^{2}}}\int_{0}^{\infty}x^{2}\widehat{g}_{1}(x^{2})dx=2C_{1}\widehat{I}=Var(Y).$$
$$\phi_{X}(y)=\frac{\rho'y}{\sqrt{2C_{1}\widehat{I}}},~\phi_{Y}(x)=\frac{\rho'x}{\sqrt{2C_{1}\widehat{I}}}.$$
\begin{equation}
\begin{split}
\widehat{I}=&\int_{0}^{\infty}t^{2}\widehat{g}_{1}(t^{2})dt=\int_{0}^{\infty}\sum_{j=0}^{N-1}
\frac{(N-1)_{j}}{j!}\frac{\Gamma(\frac{1}{2}+j)}{b^{\frac{1}{2}+j}}t^{2(N+1-j)}e^{-at^{2}}
\Phi_{2r}^{*}(-e^{-bt^{2}},\frac{1}{2},\frac{a}{b})dt\nonumber\\
=&\sum_{j=0}^{N-1}\frac{(N-1)_{j}}{j!}\frac{\Gamma(\frac{1}{2}+j)}{b^{\frac{1}{2}+j}}
\int_{0}^{\infty}t^{2(N+1-j)}e^{-at^{2}}\frac{1}{\Gamma(2r)}
\sum_{k=0}^{\infty}\frac{\Gamma(2r+k)(-1)^{k}}{k!}\frac{1}{(k+\frac{a}{b})^{\frac{1}{2}+j}}e^{-bkt^{2}}dt\\
&(setting~y=t^{2})\\
=&\sum_{j=0}^{N-1}\sum_{k=0}^{\infty}\frac{1}{2}\beta_{j}\alpha_{k}\frac{1}{(k+\frac{a}{b})^{\frac{1}{2}+j}}
\int_{0}^{\infty}y^{N-j-\frac{1}{2}}e^{-ay-bky}dy\\
=&\sum_{j=0}^{N-1}\sum_{k=0}^{\infty}\frac{1}{2}\beta_{j}\alpha_{k}\frac{1}{(k+\frac{a}{b})^{\frac{1}{2}+j}}
\frac{\Gamma(N-j+\frac{1}{2})}{(a+bk)^{N-j+\frac{1}{2}}}\\
=&\sum_{j=0}^{N-1}\sum_{k=0}^{\infty}\beta_{j}^{*}\alpha_{k}^{*},
\end{split}
\end{equation}
where $$\alpha_{k}=\frac{\Gamma(2r+k)(-1)^{k}}{k!},~\alpha_{k}^{*}=\frac{\alpha_{k}}{2(k+\frac{a}{b})^{N+1}},~
\beta_{j}=\frac{1}{\Gamma(2r)}\frac{(N-1)_{j}}{j!}\frac{\Gamma(\frac{1}{2}+j)}{b^{\frac{1}{2}+j}},~
\beta_{j}^{*}=\frac{\beta_{j}\Gamma(N-j+\frac{1}{2})}{b^{N-j+\frac{1}{2}}}.$$
The local dependence function for the elliptically symmetric generalized logistic distribution can be expressed as follows,
$$H(x,y)=\frac{\frac{N}{2b}\frac{\Phi_{2r}^{*}(-1,N+1,\frac{a}{b})}{\Phi_{2r}^{*}(-1,N,\frac{a}{b})}\rho'
+\rho'^{2}xy}{\sqrt{2C_{1}\widehat{I}+\rho'^{2}y^{2}}\sqrt{2C_{1}\widehat{I}+\rho'^{2}y^{2}}},$$
$$H(x,y)=\frac{
\frac{\frac{N}{2b}\frac{\Phi_{2r}^{*}(-1,N+1,\frac{a}{b})}{\Phi_{2r}^{*}(-1,N,\frac{a}{b})}}{2C_{1}\widehat{I}}
        \rho'+\frac{\rho'^{2}xy}{4C_{1}\widehat{I}}}
{\sqrt{1+\frac{\rho'^{2}y^{2}}{2C_{1}\widehat{I}}}\sqrt{1+\frac{\rho'^{2}y^{2}}{2C_{1}\widehat{I}}}},$$
where $$\widehat{I}=\sum_{j=0}^{\infty}\sum_{k=0}^{\infty}\beta_{j}^{*}\alpha_{k}^{*}, ~ \alpha_{k}^{*}=\frac{(-1)^{k}\Gamma(2r+k)}{2k!(k+\frac{a}{b})^{N+1}},~
~\beta_{j}^{*}=\frac{(N-1)_{j}}{j!}\frac{\Gamma(\frac{1}{2}+j)\Gamma(N-j+\frac{1}{2})}{\Gamma(2r)b^{N+1}}.$$
\section{Data analysis}
We provide a numerical illustration for the GL distribution, using data in Table 1 which is concluded by Gupta and Kundu (2010). It represents the strength measured in GPA, for single carbon fibers and impregnated 1000-carbon fiber tows. Table 1 shows the single fibers data set of 10 mm in gauge lengths with sample size 63.

If a random variable X follows a general logistic distribution then its pdf defined as
$$f(x;\theta)=\frac{e^{x/\theta}}{\theta(1+e^{x/\theta})^2},~x\in\mathbb{R}.$$
According to Chakraborty et al. (2012), the pdf of the new skew logistic (NSL) distribution is given by
$$f_{SL}(x;\lambda,\alpha,\beta)=\frac{[1+{\sin(\lambda x/(2\beta))}/\alpha]e^{-x/\beta}}{\beta[1+e^{-x/\beta}]^{2}},
~-\infty<x<\infty,~\alpha\geq1,~\lambda\in\mathbb{R},~\beta>0.$$

The GL distribution whose d.g. is defined as (2.9), is fitted to the data set and the result is compared with those for the general logistic distribution, the NSL distribution, the skew logistic (SL) distribution, the proportional reversed hazard logistic (PRHL) distribution and the exponentiated-exponential logistic (EEL) distribution. The maximum likelihood estimates, the log-likelihood value, the Akaike information criterion (AIC), the K-S test statistic and its p-value for the fitted distributions are presented in Table 2.  The results of the general logistic distribution, the NSL distribution, the SL distribution, the PRHL distribution and the EEL distribution are analyzed by Indranil and Ayman (2018). Since the data set in Table 1 is widely used for general and generalizations of logistic distributions, we consider the GL distribution with fixed $N=1.0000$, $a=1.0000$, $s=1.0000$.
The results show that the GL distribution with fixed value of $N$, $a$ and $s$ fit data better among  provided distributions in terms of Akaike information criterion. However, as for the K-S test statistic, it doesn't perform well as known distributions.

\begin{table}
 \caption{The strength in GPA for single-carbon fibers data}
 \center
 \begin{tabular}{|c|c|c|c|c|c|c|c|c|c|}
  \hline
  1.901&2.132&2.203&2.228&2.257&2.350&2.361&2.396&2.397&2.445\\\hline
  2.454&2.474&2.518&2.522&2.525&2.532&2.575&2.614&2.616&2.618\\\hline
  2.624&2.659&2.675&2.738&2.740&2.856&2.917&2.928&2.937&2.937\\\hline
  2.977&2.996&3.030&3.125&3.139&3.145&3.220&3.223&3.235&3.243\\\hline
  3.264&3.272&3.294&3.332&3.346&3.377&3.408&3.435&3.493&3.501\\\hline
  3.537&3.554&3.562&3.628&3.852&3.871&3.886&3.971&4.024&4.027\\\hline
  4.225&4.395&5.020&     &     &     &     &     &     &     \\
  \hline
\end{tabular}
\end{table}

\begin{table}
 \center
 \caption{Parameter estimates}
\resizebox{\textwidth}{!}{\begin{tabular}{|c|c|c|c|c|c|l|}
  \hline
  Distribution & Logistic & NSL & SL & PRHL & EEL & GL \\ \hline
  Parameter estimates & $\widehat{\theta}=0.19975$ & $\widehat{\alpha}=1.0000$ & $\widehat{\alpha}=0.5400$ & $\widehat{\alpha}=3.2761$ & $\widehat{\alpha}=218.2300$ & $\widehat{N}=1.0000$\\
                      &                            & $\widehat{\beta}=1.9975$  & $\widehat{\lambda}=2.6800$  & $\widehat{\lambda}=2.2192$ & $\widehat{\lambda}=0.0946$ & $\widehat{a}=1.0000$  \\
                      &                            & $\widehat{\lambda}=1.9694$ & $\widehat{\mu}=2.740$ & $\widehat{\mu}=2.3369$ & $\widehat{\theta}=0.0486$& $\widehat{s}=1.0000$\\
                      &                            &                            &   &   &   & $\widehat{b}=8.7827e+04$\\
                      &                            &                            &   &   &   & $\widehat{\mu}=3.0593$\\
                      &                            &                            &   &   &   & $\widehat{\sigma^2}=0.7588$\\
                      &                            &                            &   &   &   & $\widehat{r}=4.1739e-38$\\\hline
  Log likelihood & -165.5826 & -123.4458 & -58.0299 & -58.9896 & -56.8643 &  -49.6587 \\\hline
  AIC& 333.1652 & 248.8916 & 122.0597 & 119.9792 & 119.7286 &  107.3174 \\\hline
  K-S& 0.7123 & 0.5532 & 0.0918 & 0.0844 & 0.0735 & 0.0987 \\\hline
  K-S p-value& 0.0000 & 0.0000 & 0.6632 & 0.7603 & 0.8853 &  0.5714\\
  \hline
\end{tabular}}
\end{table}

\clearpage
\section{Concluding remarks}
This paper defined the generalized logistic distribution whose density generator is defined as
$$ g(t)=\frac{t^{N-1}\exp(-at^{s_{1}})}{\left(1+\exp(-bt^{s_{2}})\right)^{2r}},~t>0,\eqno(9.1)$$ where $2N+n>2$, $a,~b,~ s_{1},~ s_{2}>0$, $r\geq0$ are constants. By setting different $a,~b,~s_{1},~s_{2},~r,~N$ in $(9.1)$, we obtained various density generators of elliptical distributions, such as the normal distribution, the Kotz type distribution, the  exponential power distribution, the symmetric logistic distribution and generalized logistic type \RNum{1}, \RNum{3}, \RNum{4} distribution, etc. Our interest is to study the inconsistency properties and various probabilistic properties of this distribution including marginal distributions, conditional distributions, linear transformations, characteristic functions. In addition, we gave a data analysis which shows that the GL distributions are more flexible than other distributions. We would give further research on statistic inference of this new kind of elliptical distributions in the subsequent research.
\newpage
\section*{Appendix A. Proofs}
\subsection*{Appendix A.1. Proof of $(5.1)-(5.4)$}
\textbf{Proof.} When the dimension $n>1$, $\Omega_{n}(\|t\|^{2})$ defined as $$\Omega_{n}(\|t\|^{2})=\frac{\Gamma(\frac{n}{2})}{\sqrt{\pi}}\sum_{k=0}^{\infty}\frac{(-1)^{k}\|t\|^{2k}}{(2k)!}
   \frac{\Gamma(\frac{2k+1}{2})}{\Gamma(\frac{n+2k}{2})},$$ we have the characteristic function as follows.
\begin{equation}
\begin{split}
\psi_{\textbf{X}}(\textbf{t})=&E(e^{i\textbf{t}^{\emph{T}}\textbf{X}})
=e^{i\textbf{t}^{\emph{T}}\bm{\mu}}E(e^{i\textbf{t}^{\emph{T}}R\mathbf{\Sigma}^{\frac{1}{2}}\textbf{U}^{(n)}})
=e^{i\textbf{t}^{\emph{T}}\bm{\mu}}E[\Omega_{n}(R^{2}\textbf{t}^{\emph{T}}\mathbf{\Sigma}\textbf{t})]\nonumber\\
=&e^{i\textbf{t}^{\emph{T}}\bm{\mu}}\int_{0}^{\infty}\Omega_{n}(v\textbf{t}^{\emph{T}}\mathbf{\Sigma}\textbf{t})
\frac{1}{\int_{0}^{\infty}t^{\frac{n}{2}-1}g_{n}(t)dt}v^{\frac{n}{2}-1}\frac{v^{N-1}\exp(-av^{s})}{(1+\exp(-bv^{s}))^{2r}}dv\\
=&e^{i\textbf{t}^{\emph{T}}\bm{\mu}}\left[\int_{0}^{\infty}\frac{t^{\frac{n}{2}+N-2}\exp(-at^{s})}{(1+\exp(-bt^{s}))^{2r}}\right]^{-1} \int_{0}^{\infty}\Omega_{n}(v\textbf{t}^{\emph{T}}\mathbf{\Sigma}\textbf{t})\frac{v^{\frac{n}{2}+N-2}\exp(-av^{s})}{(1+\exp(-bv^{s}))^{2r}}dv\\
=&e^{i\textbf{t}^{\emph{T}}\bm{\mu}}\left[\frac{\Gamma(\frac{1}{s}(\frac{n}{2}+N-1))}{b^{\frac{1}{s}(\frac{n}{2}+N-1)}s}
\Phi_{2r}^{*}(-1,\frac{1}{s}(\frac{n}{2}+N-1),\frac{a}{b})\right]^{-1}\\
&\times\int_{0}^{\infty}\frac{\Gamma(\frac{n}{2})}{\pi^{\frac{1}{2}}}\sum_{k=0}^{\infty}
\frac{(-1)^{k}(\textbf{t}^{\emph{T}}\mathbf{\Sigma}\textbf{t})^{k}v^{k}\Gamma(\frac{2k+1}{2})}{(2k)!\Gamma(\frac{n+2k}{2})}
\frac{v^{\frac{n}{2}+N-2}\exp(-av^{s})}{(1+\exp(-bv^{s}))^{2r}}dv\\
=&e^{i\textbf{t}^{\emph{T}}\bm{\mu}}\left[\frac{\Gamma(\frac{1}{s}(\frac{n}{2}+N-1))}{b^{\frac{1}{s}(\frac{n}{2}+N-1)}s}
\Phi_{2r}^{*}(-1,\frac{1}{s}(\frac{n}{2}+N-1),\frac{a}{b})\right]^{-1}\\
 &\times\sum_{k=0}^{\infty}\frac{\Gamma(\frac{n}{2})}{\pi^{\frac{1}{2}}}
 \frac{(-1)^{k}(\textbf{t}^{\emph{T}}\mathbf{\Sigma}\textbf{t})^{k}\Gamma(\frac{2k+1}{2})}{(2k)!\Gamma(\frac{n+2k}{2})}
 \int_{0}^{\infty}\frac{v^{\frac{n}{2}+N+k-2}\exp(-av^{s})}{(1+\exp(-bv^{s}))^{2r}}dv\\
=&e^{it^{\emph{T}}\bm{\mu}}\left[\frac{\Gamma(\frac{1}{s}(\frac{n}{2}+N-1))}{b^{\frac{1}{s}(\frac{n}{2}+N-1)}s}
\Phi_{2r}^{*}(-1,\frac{1}{s}(\frac{n}{2}+N-1),\frac{a}{b})\right]^{-1}\\
&\times\sum_{k=0}^{\infty}\frac{\Gamma(\frac{n}{2})}{\pi^{\frac{1}{2}}}
 \frac{(-1)^{k}(\textbf{t}^{\emph{T}}\mathbf{\Sigma}\textbf{t})^{k}\Gamma(\frac{2k+1}{2})}{(2k)!\Gamma(\frac{n+2k}{2})}
 \frac{\Gamma(\frac{1}{s}(\frac{n}{2}+N-1))}{b^{\frac{1}{s}(\frac{n}{2}+N-1)}s}
\Phi_{2r}^{*}(-1,\frac{1}{s}(\frac{n}{2}+N-1),\frac{a}{b})\\
=&e^{i\textbf{t}^{\emph{T}}\bm{\mu}}\frac{\Gamma(\frac{n}{2})}{\pi^{\frac{1}{2}}}\sum_{k=0}^{\infty}
 \frac{(-1)^{k}}{(2k)!}\frac{\Gamma(k+\frac{1}{2})\Gamma(\frac{1}{s}(\frac{n}{2}+N+k-1))}{\Gamma(k+\frac{n}{2})\Gamma(\frac{1}{s}(\frac{n}{2}+N-1))b^{\frac{1}{s}(\frac{n}{2}+N-1)}}\\
 &\times\frac{\Phi_{2r}^{*}(-1,\frac{1}{s}(\frac{n}{2}+N+k-1),\frac{a}{b})}{\Phi_{2r}^{*}(-1,\frac{1}{s}(\frac{n}{2}+N-1),\frac{a}{b})}
 (\textbf{t}^{\emph{T}}\mathbf{\Sigma}\textbf{t})^{k}\\
 =&e^{i\textbf{t}^{\emph{T}}\bm{\mu}}\sum_{k=0}^{\infty}\frac{(-1)^{k}}{(2k)!}\gamma_{k}
 \frac{\Phi_{2r}^{*}(-1,\frac{1}{s}(\frac{n}{2}+N+k-1),\frac{a}{b})}{\Phi_{2r}^{*}(-1,\frac{1}{s}(\frac{n}{2}+N-1),\frac{a}{b})}
 (\textbf{t}^{\emph{T}}\mathbf{\Sigma}\textbf{t})^{k},
 \end{split}
\end{equation}
where$$\gamma_{k}(N,n,b,s)=\frac{\Gamma(\frac{n}{2})}{\pi^{\frac{1}{2}}b^{\frac{1}{s}(\frac{n}{2}+N-1)}}
\frac{\Gamma(k+\frac{1}{2})\Gamma(\frac{1}{s}(\frac{n}{2}+N+k-1))}{\Gamma(k+\frac{n}{2})\Gamma(\frac{1}{s}(\frac{n}{2}+N-1))}.$$
Here $\Phi_{2r}^{*}$ is the generalized Hurwitz-Lerch zeta function.

When the dimension $n>1$, $\Omega_{n}(\|t\|^{2})$ defined as $$\Omega_{n}(\|t\|^{2})=_{0}F_{1}(\frac{n}{2};-\frac{1}{4}\|t\|^{2}),$$ we have the characteristic function as follows.
\begin{equation}
\begin{split}
\psi_{\textbf{X}}(\textbf{t})=&e^{i\textbf{t}^{\emph{T}}\bm{\mu}}E[\Omega_{n}(R^{2}\textbf{t}'\mathbf{\Sigma}\textbf{t})]\nonumber\\
 =&e^{i\textbf{t}^{\emph{T}}\bm{\mu}}\left[\int_{0}^{\infty}\frac{t^{\frac{n}{2}+N-2}\exp(-at^{s})}{(1+\exp(-bt^{s}))^{2r}}\right]^{-1} \int_{0}^{\infty}\Omega_{n}(v\textbf{t}^{\emph{T}}\mathbf{\Sigma}\textbf{t})\frac{v^{\frac{n}{2}+N-2}\exp(-av^{s})}{(1+\exp(-bv^{s}))^{2r}}dv\\
=&e^{i\textbf{t}^{\emph{T}}\bm{\mu}}\left[\frac{\Gamma(\frac{1}{s}(\frac{n}{2}+N-1))}{b^{\frac{1}{s}(\frac{n}{2}+N-1)}s}
\Phi_{2r}^{*}(-1,\frac{1}{s}(\frac{n}{2}+N-1),\frac{a}{b})\right]^{-1}\\
&\times\int_{0}^{\infty} {}_{0}F_{1}(\frac{n}{2};\frac{1}{4}v\textbf{t}^{\emph{T}}{\Sigma}\textbf{t})\frac{v^{\frac{n}{2}+N-2}\exp(-av^{s})}{(1+\exp(-bv^{s}))^{2r}}dv\\
=&e^{i\textbf{t}^{\emph{T}}\bm{\mu}}\left[\frac{\Gamma(\frac{1}{s}(\frac{n}{2}+N-1))}{b^{\frac{1}{s}(\frac{n}{2}+N-1)}s}
\Phi_{2r}^{*}(-1,\frac{1}{s}(\frac{n}{2}+N-1),\frac{a}{b})\right]^{-1}\\
&\times\sum_{k=0}^{\infty}\frac{\Gamma(\frac{1}{s}(\frac{n}{2}+N-1))}{b^{\frac{1}{s}(\frac{n}{2}+N-1)}s}
\frac{(\textbf{t}^{\emph{T}}\Sigma \textbf{t})^{k}}{4^{k}(\frac{n}{2})^{[k]}k!}\Phi_{2r}^{*}(-1,\frac{1}{s}(\frac{n}{2}+N+k-1),\frac{a}{b})\\
=&e^{i\textbf{t}^{\emph{T}}\bm{\mu}}\sum_{k=0}^{\infty}
\frac{\Gamma(\frac{1}{s}(\frac{n}{2}+N+k-1))}{\Gamma(\frac{1}{s}(\frac{n}{2}+N-1))b^{\frac{k}{s}}4^{k}(\frac{n}{2})^{[k]}k!}
\frac{\Phi_{2r}^{*}(-1,\frac{1}{s}(\frac{n}{2}+N+k-1),\frac{a}{b})}{\Phi_{2r}^{*}(-1,\frac{1}{s}(\frac{n}{2}+N-1),\frac{a}{b})}
(\textbf{t}^{\emph{T}}\mathbf{\Sigma} \textbf{t})^{k}.
\end{split}
\end{equation}
Thus, we have characteristic generators as follows:
\begin{equation}
\begin{split}
\phi_{\textbf{X}}(\|\textbf{u}_{n}\|^{2})
 =&\sum_{k=0}^{\infty}
 \frac{(-1)^{k}}{(2k)!}\gamma_{k}(N,n,b,s)
 \frac{\Phi_{2r}^{*}(-1,\frac{1}{s}(\frac{n}{2}+N+k-1),\frac{a}{b})}{\Phi_{2r}^{*}(-1,\frac{1}{s}(\frac{n}{2}+N-1),\frac{a}{b})}
 (\|\textbf{u}_{n}\|)^{2k},\nonumber\\
 \phi_{\textbf{X}}(\|\textbf{u}_{n}\|^{2})
 =&\sum_{k=0}^{\infty}
\frac{\Gamma(\frac{1}{s}(\frac{n}{2}+N+k-1))}{\Gamma(\frac{1}{s}(\frac{n}{2}+N-1))b^{\frac{k}{s}}4^{k}(\frac{n}{2})^{[k]}k!}
\frac{\Phi_{2r}^{*}(-1,\frac{1}{s}(\frac{n}{2}+N+k-1),\frac{a}{b})}{\Phi_{2r}^{*}(-1,\frac{1}{s}(\frac{n}{2}+N-1),\frac{a}{b})}
(\|\textbf{u}_{n}\|)^{2k},
 \end{split}
\end{equation}
where$$\gamma_{k}(N,n,b,s)=\frac{\Gamma(\frac{n}{2})}{\pi^{\frac{1}{2}}b^{\frac{1}{s}(\frac{n}{2}+N-1)}}
\frac{\Gamma(k+\frac{1}{2})\Gamma(\frac{1}{s}(\frac{n}{2}+N+k-1))}{\Gamma(k+\frac{n}{2})\Gamma(\frac{1}{s}(\frac{n}{2}+N-1))}.$$
\subsection*{Appendix A.2. Proof of $(7.1)-(7.4)$}
\textbf{Proof.} When the dimension $n>1$, $\Omega_{n}(\|t\|^{2})$ defined as
$$\Omega_{n}(\|t\|^{2})=\frac{\Gamma(\frac{n}{2})}{\sqrt{\pi}}\sum_{k=0}^{\infty}\frac{(-1)^{k}\|t\|^{2k}}{(2k)!}
   \frac{\Gamma(\frac{2k+1}{2})}{\Gamma(\frac{n+2k}{2})},$$ we have the characteristic function as follows.
\begin{equation}
\begin{split}
\psi_{\textbf{Y}}(\textbf{t})=&E(e^{i\textbf{t}^{\emph{T}}\textbf{Y}})
=e^{i\textbf{t}^{\emph{T}}\bm{\mu}_{\textbf{Y}}}\frac{1}{\int_{0}^{\infty}t^{\frac{m}{2}-1}g_{m}(t)dt}\\ \nonumber
  &\times\int_{0}^{\infty}\Omega_{n}(v\textbf{t}^{\emph{T}}\mathbf{\Sigma_{\textbf{Y}}}\textbf{t})
  \sum_{j=0}^{N-1}\frac{(N-1)_{j}}{j!}\frac{\Gamma(\frac{n-m}{2}+j)}{b^{\frac{n-m}{2}+j}}v^{N+\frac{m}{2}-j-2}e^{-av}
  \Phi_{2r}^{*}(-e^{-bv},\frac{n-m}{2}+j,\frac{a}{b})dv\\
=&e^{it^{\emph{T}}\bm{\mu}_{\textbf{Y}}}\left[\sum_{k=0}^{\infty}\sum_{j=0}^{N-1}\alpha_{k}^{*}\beta_{j}
\Gamma(\frac{m}{2}+N-j-1)\right]^{-1}I_{3},
\end{split}
\end{equation}
where $$\beta_{j}=\frac{(N-1)_{j}}{j!}\frac{\Gamma(\frac{n-m}{2}+j)}{b^{\frac{n-m}{2}+j}\Gamma(2r)},~
\alpha_{k}=\frac{\Gamma(2r+k)(-1)^{k}}{k!},~
\alpha_{k}^{*}=\frac{\alpha_{k}}{(k+\frac{a}{b})^{\frac{m}{2}+N-1}}.
$$
\begin{equation}
\begin{split}
I_{3}=&\int_{0}^{\infty}\frac{\Gamma(\frac{m}{2})}{\pi^{\frac{1}{2}}}\sum_{l=0}^{\infty}
 \frac{(-1)^{l}(\textbf{t}^{\emph{T}}\mathbf{\Sigma_{\textbf{Y}}}\textbf{t})^{l}v^{l}\Gamma(\frac{2l+1}{2})}{(2l)!\Gamma(\frac{m+2l}{2})}\nonumber\\
&\times\sum_{k=0}^{N-1}\frac{(N-1)_{k}}{k!}\frac{\Gamma(\frac{n-m}{2}+k)}{b^{\frac{n-m}{2}+k}}v^{N+\frac{m}{2}-k-2}e^{-av}
  \Phi_{2r}^{*}(-e^{-bv},\frac{n-m}{2}+k,\frac{a}{b})dv\\
=&\frac{\Gamma(\frac{m}{2})}{\pi^{\frac{1}{2}}}\sum_{l=0}^{\infty}\frac{(-1)^{l}\Gamma(\frac{2l+1}{2})}{(2l)!\Gamma(\frac{m+2l}{2})}
(\textbf{t}^{\emph{T}}\mathbf{\Sigma_{\textbf{Y}}}\textbf{t})^{l}\sum_{k=0}^{N-1}\beta_{k}\int_{0}^{\infty}v^{N+\frac{m}{2}-k-2+l}e^{-av}
 \sum_{p=0}^{\infty}\frac{\Gamma(2r+p)}{p!}\frac{(-1)^{p}e^{-bpv}}{(p+\frac{a}{b})^{\frac{n-m}{2}+k}}dv\\
=&\sum_{l=0}^{\infty}\sum_{k=0}^{N-1} \sum_{p=0}^{\infty}
\frac{(-1)^{l}}{(2l)!}\frac{\Gamma(\frac{m}{2})\Gamma(\frac{2l+1}{2})}{\Gamma(\frac{1}{2})\Gamma(\frac{m+2l}{2})}\beta_{k}\alpha_{p}
   \frac{\Gamma(\frac{m}{2}+N-k+l-1)}{b^{N-k+l}(p+\frac{a}{b})^{\frac{n-m}{2}+N+l}}(\textbf{t}^{\emph{T}}\mathbf{\Sigma_{\textbf{Y}}}\textbf{t})^{l}\\
=&\sum_{l=0}^{\infty}\sum_{k=0}^{N-1} \sum_{p=0}^{\infty}
\frac{(-1)^{l}}{(2l)!}\alpha_{p}^{*}\beta_{k}A_{l}\Gamma(\frac{m}{2}+N-k+l-1)(\textbf{t}^{\emph{T}}\mathbf{\Sigma_{\textbf{Y}}}\textbf{t})^{l},
\end{split}
\end{equation}
\begin{equation}
\begin{split}
\psi_{\textbf{Y}}(\textbf{t})=&e^{i\textbf{t}^{\emph{T}}\bm{\mu}_{\textbf{Y}}}
\frac{\sum_{l=0}^{\infty}\sum_{k=0}^{N-1} \sum_{p=0}^{\infty}
\frac{(-1)^{l}}{(2l)!}\alpha_{p}^{*}\beta_{k}A_{l}\Gamma(\frac{m}{2}+N-k+l-1)}
{\sum_{k=0}^{\infty}\sum_{j=0}^{N-1}\alpha_{k}^{*}\beta_{j}\Gamma(\frac{m}{2}+N-j-1)}(\textbf{t}^{\emph{T}}\mathbf{\Sigma_{\textbf{Y}}}\textbf{t})^{l},\nonumber\\
 \phi_{(m),y}(\|\bm{\xi}_{(m)}\|^{2})
=&\frac{\sum_{l=0}^{\infty}\sum_{k=0}^{N-1} \sum_{p=0}^{\infty}
\frac{(-1)^{l}}{(2l)!}\alpha_{p}^{*}\beta_{k}A_{l}\Gamma(\frac{m}{2}+N-k+l-1)}
{\sum_{k=0}^{\infty}\sum_{j=0}^{N-1}\alpha_{k}^{*}\beta_{j}\Gamma(\frac{m}{2}+N-j-1)}\|\bm{\xi}_{(m)}\|^{2l},
 \end{split}
\end{equation}
where $$A_{x}\triangleq A_{x}(N,m,a,b,p,k)=\frac{B(\frac{m}{2},x+\frac{1}{2})}{B(\frac{m}{2}+x,\frac{1}{2})b^{N-k+x}(p+\frac{a}{b})^{x-\frac{m}{2}+1}},$$
$$\beta_{x}\triangleq\beta_{x}(N,n,m,b,r)=\frac{(N-1)_{x}}{x!}\frac{\Gamma(\frac{n-m}{2}+x)}{b^{\frac{n-m}{2}+x}\Gamma(2r)},$$
$$\alpha_{x}(r)\triangleq\alpha_{x}(r)=\frac{\Gamma(2r+x)(-1)^{x}}{x!},~\alpha_{x}^{*}\triangleq\alpha_{x}^{*}(N,n,a,b,r)=\frac{\alpha_{x}}{(x+\frac{a}{b})^{\frac{n}{2}+N-1}}.
$$
When the dimension $n>1$, $\Omega_{n}(\|t\|^{2})$ defined as
$$\Omega_{n}(\|t\|^{2})=_{0}F_{1}(\frac{n}{2};-\frac{1}{4}\|t\|^{2}),$$
we have the characteristic function as follows.
\begin{equation}
\begin{split}
\psi_{\textbf{Y}}(\textbf{t})=&e^{i\textbf{t}^{\emph{T}}\bm{\mu}_{\textbf{Y}}}\frac{1}{\int_{0}^{\infty}t^{\frac{m}{2}-1}g_{m}(t)dt} \int_{0}^{\infty}\Omega_{n}(v\textbf{t}^{\emph{T}}\mathbf{\Sigma_{\textbf{Y}}}\textbf{t})\nonumber\\
  &\times\sum_{j=0}^{N-1}\frac{(N-1)_{j}}{j!}\frac{\Gamma(\frac{n-m}{2}+j)}{b^{\frac{n-m}{2}+j}}v^{N+\frac{m}{2}-j-2}e^{-av}
  \Phi_{2r}^{*}(-e^{-bv},\frac{n-m}{2}+j,\frac{a}{b})dv\\
=&e^{i\textbf{t}^{\emph{T}}\bm{\mu}_{\textbf{Y}}}\left[\sum_{k=0}^{\infty}\sum_{j=0}^{N-1}\alpha_{k}^{*}\beta_{j}
\Gamma(\frac{m}{2}+N-j-1)\right]^{-1}\int_{0}^{\infty}\frac{\Gamma(\frac{m}{2})}{\pi^{\frac{1}{2}}}
\sum_{l=0}^{\infty}\frac{(\textbf{t}^{\emph{T}}\mathbf{\Sigma}_{\textbf{Y}}\textbf{t})^{l}}{4^{l}(\frac{m}{2})^{[l]}}\\
&\times\sum_{j=0}^{N-1}\frac{(N-1)_{j}}{j!}\frac{\Gamma(\frac{n-m}{2}+j)}{b^{\frac{n-m}{2}+j}}v^{N+\frac{m}{2}-j-2}e^{-av}
 \Phi_{2r}^{*}(-e^{-bv},\frac{n-m}{2}+j,\frac{a}{b})dv\\
=&e^{i\textbf{t}^{\emph{T}}\bm{\mu}_{\textbf{Y}}}\frac{\sum_{l=0}^{\infty}\sum_{j=0}^{N-1}\frac{\Gamma(\frac{m}{2})}{\pi^{\frac{1}{2}}}
\frac{(\textbf{t}^{\emph{T}}\mathbf{\Sigma}_{\textbf{Y}}\textbf{t})^{l}}{4^{l}(\frac{m}{2})^{[l]}}
\sum_{k=0}^{\infty}\beta_{j}\alpha_{k}\frac{\Gamma(N+\frac{m}{2}-j-1)}{b^{N+\frac{m}{2}-j-1}} \Phi_{2r}^{*}(-e^{-bv},\frac{n-m}{2}+j,\frac{a}{b})}{\sum_{k=0}^{\infty}\sum_{j=0}^{N-1}\alpha_{k}^{*}\beta_{j}
\Gamma(\frac{m}{2}+N-j-1)}\nonumber\\
=&e^{i\textbf{t}^{\emph{T}}\bm{\mu}_{\textbf{Y}}}\frac{\sum_{l=0}^{\infty}\sum_{j=0}^{N-1}\sum_{k=0}^{\infty}\alpha_{k}\beta_{j}^{*}
\Phi_{2r}^{*}(-e^{-bv},\frac{n-m}{2}+j,\frac{a}{b})}{\sum_{k=0}^{\infty}\sum_{j=0}^{N-1}\alpha_{k}^{*}\beta_{j}\Gamma(\frac{m}{2}+N-j-1)}
\frac{\Gamma(\frac{m}{2})}{\pi^{\frac{1}{2}}4^{l}(\frac{m}{2})^{[l]}}(\textbf{t}^{\emph{T}}\mathbf{\Sigma}_{\textbf{Y}}\textbf{t})^{l},
\end{split}
\end{equation}
$$
\phi_{(m),y}(\|\bm{\xi}_{(m)}\|^{2})=\frac{\sum_{l=0}^{\infty}\sum_{j=0}^{N-1}\sum_{k=0}^{\infty}\alpha_{k}\beta_{j}^{*}
\Phi_{2r}^{*}(-e^{-bv},\frac{n-m}{2}+j,\frac{a}{b})}{\sum_{k=0}^{\infty}\sum_{j=0}^{N-1}\alpha_{k}^{*}\beta_{j}\Gamma(\frac{m}{2}+N-j-1)}
\frac{\Gamma(\frac{m}{2})}{\pi^{\frac{1}{2}}4^{l}(\frac{m}{2})^{[l]}}\|\bm{\xi}_{(m)}\|^{2l},
$$
 where
 $$\beta_{x}\triangleq\beta_{x}(N,n,m,r)=\frac{(N-1)_{x}}{x!}\frac{\Gamma(\frac{n-m}{2}+x)}{\Gamma(2r)b^{\frac{n-m}{2}+x}},~
  \beta_{x}^{*}\triangleq\beta_{x}^{*}(N,n,m,b,r)=\frac{\beta_{x}\Gamma(N+\frac{m}{2}-x-1)}{b^{N+\frac{m}{2}-x-1}},$$
  $$\alpha_{x}\triangleq\alpha_{x}(r)=\frac{(-1)^{x}\Gamma(2r+x)}{x!},~
  \alpha_{x}^{*}\triangleq\alpha_{x}^{*}(N,m,a,b,r)=\frac{\alpha_{x}}{(x+\frac{a}{b})^{\frac{n}{2}+N-1}}.$$
Here $\Phi_{2r}^{*}$ is the generalized Hurwitz-Lerch zeta function.
\section*{Acknowledgements}
The authors are grateful to the referee and the Editor whose comments and suggestions greatly improved
the article.

\noindent The research was supported by the National Natural Science Foundation of China (No. 11171179, 11571198).


\newpage
\begin{thebibliography}{11}
\bibitem{1} Ali MM, Mikhail NN, Haq MS (1978) A class of bivariate distributions including the bivariate logistic.  Journal of Multivariate Analysis 8:405-412.
\bibitem{2} Arashi M, Nadarajh S (2016) Generalized elliptical distributions. Communications in Statistics-Theory and Methods 46:6412-6432.
\bibitem{3}Bairmov I, Kotz S, Kozubowski TJ (2003) A new measure of liner local dependence. Statistics 37(3):243-258.
\bibitem{4} Balakrishnan N, Leung MY (2007) Order statistics from the type $\RNum{1}$ generalized logistic distribution. Communications in Statistics-Simulation and Computation 17(1):25-50.
\bibitem{5} Chakraborty S, Hazarika PJ, Ali MM (2012) A new skew logistic distribution and its properties. Pakistan Journal of Statistics 28(4):513-524.
\bibitem{6} Cambanis S, Huang S, Simons G (1981) On the theory of elliptically contoured distributions. Journal of Multivaruate Analysis 11:368-385.
\bibitem{7} Denuit M, Dhaene J, Goovaerts M, Kaas R (2005) Actuarial Theory for Dependent Risks Measure, Orders and Models. John Wiley Sons, Ltd.
\bibitem{8} Fang KT, Kotz S, Ng KW (1990) Symmetric Multivariate and Related Distributions. Chapman and Hall.
\bibitem{9} Ghosh I, Alzaatreh A (2018) A new class of generalized logistic distribution. Communications in Statistics - Theory and Methods 47(9):2043-2055.
\bibitem{10}Gupta RD, Kundu D (2010) Generalized logistic distributions. Journal of Applied Statistical Science 18(1):1-23.
\bibitem{11}Jensen DR (1985) Multivariate distributions. Encyclopedia of Statistical Sciences 6 (eds S. Kotz, N. L. Johnson, and C. B. Read):43-45.
\bibitem{12} Kano Y (1994) Consistency property of elliptical probability density functions. Journal of Multivaruate Analysis 51:139-147.
\bibitem{13}  Kotz S, Ostrovskii I (1994) Characteristic functions of a class of elliptical distributions. Journal of Multivariate Analysis 49:164-178.
\bibitem{14} Kotz S (1975) Multivariate distributions at a cross road. Stastistical Distributions in Scientific Work 1(eds G. P. Patil et al.):247-270.
\bibitem{15} Kotz S, Nadarajah S (2001) Some extremal type elliptical distributions. Statistics $\&$ Probability Letters 54:171-182.
\bibitem{16} Kotz S, Nadarajah S (2003) Local dependence functions for the elliptically symmetric distributions. $Sankhy\bar{a}$: The Indian Journal of Statistics 65:207-223.
\bibitem{17} Landsman Z, Valdez EA (2003) Tail conditional expectations for elliptical distributions. North American Actuarial Journal 7(4):55-71.
\bibitem{18} Landsman Z, Makov U, Shushi T (2016a) Multivariate tail conditional expectation for elliptical distributions. Insurance: Mathematics and Economics 70:216-223.
\bibitem{19} Landsman Z, Makov U, Shushi T (2016b) Tail conditional moments for elliptical and log-elliptical distributions and applications. Insurance: Mathematics and Economics 71:179-188.
\bibitem{20} Landsman Z, Makov U, Shushi T (2018) A multivariate tail covariance measure for elliptical distributions. Insurance: Mathematics and Economics 81:27-35.
\bibitem{21} Landsman Z, Vanduffel S, Yao J (2014) Some Stein-type inequalities for multivariate elliptical distributions and applications. Statistics $\&$ Probability Letters 97:54-62.
\bibitem{22}  Liang JJ, Bentler PM (1998) Characterizations of some subclasses of spherical distributions. Statistics $\&$ Probability Letters 40:155-164.
\bibitem{23} Lin SD, Srivastava HM, Wang PY (2006) Some expansion formulas for a class of generalized Hurwitz-Lerch Zeta functions. Integral Transforms and Special Functions 17:817-827.
\bibitem{24}  Nadarajah S (2003) The Kotz-type distribution with applications. Statistics: A Journal of Theoretical and Applied Statistics 37(4):341-358.
\bibitem{25} Nassar MM, Elmasry A (2012) A study of generalized logistic distribution. Journal of the Egyptian Mathematical Society 20:126-133.
\bibitem{26} Nadarajah S, Kotz S (2005) A generalized logistic distribution. International Journal of Mathematics and Mathematical Sciences 9:3169-3174.
\bibitem{27} Nelson RB (2006) An Introduction to Copulas. Springer Science+Business Media, Inc.
\bibitem{28}Rahman G, Nisar KS, Arshad MA. A new extension of Hurwitz-Lerch Zeta function, arXiv:1802.07823vl.
\bibitem{29} Schmidt R (2002) Tail dependence for elliptically contoured distributions. Mathematical Methods of Operations Research 55:301-327.
\bibitem{30}Yeh HC (2009) Multivariate semi-logistic distributions. Journal of Multivariate Analysis 101:893-908.
\bibitem{31}Yin CC, Sha XY (2019) A new class of symmetric distributions including the elliptically symmetric logistic, arXiv:1810.10692.

\end{thebibliography}
\end{document}